\documentclass[12pt]{amsart}
\usepackage[T1]{fontenc}
\usepackage{lmodern,amsmath,amssymb,amsthm,mathtools,microtype}

\usepackage[margin=4cm]{geometry}

\usepackage{latexsym}
\usepackage{epsfig}
\usepackage{amsfonts}
\usepackage{enumerate}
\usepackage{times,float}
\usepackage{graphics}

\usepackage[dvipsnames]{xcolor}

\usepackage[nodisplayskipstretch]{setspace}
\usepackage{enumitem,needspace}
\usepackage[dvipsnames]{xcolor}
\usepackage[pagebackref=true,colorlinks=true,linkcolor=MidnightBlue,citecolor=MidnightBlue,urlcolor=MidnightBlue]{hyperref}

\newtheorem{theorem}{Theorem}[section]
\newtheorem{lemma}[theorem]{Lemma}
\newtheorem{proposition}[theorem]{Proposition}
\newtheorem{corollary}[theorem]{Corollary}
\theoremstyle{definition}
\newtheorem{definition}[theorem]{Definition}
\theoremstyle{remark}
\newtheorem{remark}[theorem]{Remark}
\newtheorem*{example*}{Example}
\newcommand{\A}{\mathbb A}
\newcommand{\PP}{\mathbb P}
\newcommand{\C}{\mathbb C}
\newcommand{\Q}{\mathbb Q}

\newcommand{\Gm}{\mathbb G_m}
\newcommand{\cO}{\mathcal O}
\newcommand{\OO}{\mathcal O}
\newcommand{\HH}{\mathbb H}
\newcommand{\De}{D_{\epsilon}}
\newcommand{\Tot}{\operatorname{Tot}}
\newcommand{\Proj}{\operatorname{Proj}}
\newcommand{\Spec}{\operatorname{Spec}}
\newcommand{\Sym}{\operatorname{Sym}}
\newcommand{\rk}{\operatorname{rk}}
\newcommand{\Isom}{\operatorname{Isom}}
\newcommand{\Supp}{\operatorname{Supp}}

\newcommand{\red}{\mathrm{red}}
\newcommand{\id}{\mathrm{id}}
\newcommand{\ord}{\operatorname{ord}}
\newcommand{\Hom}{\operatorname{Hom}}
\newcommand{\Nef}{\operatorname{Nef}}
\newcommand{\Int}{\operatorname{Int}}

\newcommand{\Exc}{\operatorname{Exc}}
\newcommand{\lcm}{\operatorname{lcm}}
\newcommand{\PGL}{\operatorname{PGL}}

\numberwithin{equation}{section}

\setlist[enumerate]{label=\textup{(\roman*)},leftmargin=*,itemsep=3pt,topsep=5pt}

\title{On rational curves relative to a divisor}
\author{Jiun-Cheng Chen}
\address{Department of Mathematics\\ Third General Building\\ National Tsing-Hua University\\ No. 101 Sec. 2 Kuang Fu Road\\ Hsinchu, 30043\\ Taiwan}
\email{jcchen@math.nthu.edu.tw}

\author[Hsian-Hua Tseng]{Hsian-Hua Tseng}
\address{Department of Mathematics\\ Ohio State University\\ 100 Math Tower, 231 West 18th Ave.\\ Columbus, OH 43210\\ USA}
\email{hhtseng@math.ohio-state.edu}
\date{\today}

\begin{document}
\begin{abstract}
Let $X$ be a smooth complex projective variety and $D$ a reduced simple normal crossings divisor. We prove that, given a complete curve of negative $(K_X+D)$-degree, through its general points there is a rational curve whose normalization meets $D$ in at most one point. The proof uses ordinary weighted jets of a smooth family constructed from $(X,D)$, a positive parameter line bundle, and an intersection calculation on a ruled surface. We obtain an explicit degree bound depending on a projective embedding.

We apply this construction to a cyclic cover compactification associated with a smooth Campana orbifold $(X,\sum_i(1-1/m_i)D_i)$. If its anticanonical divisor is big, we show that there is a covering family of free divisible Campana rational curves with fixed numerical class and contact data. In Picard number one case, the orbifold is Campana rationally connected. 
\end{abstract}
\maketitle
\vspace{-1.2cm}
\tableofcontents

\section{Introduction}\label{sec:introduction}
We work over $\C$. 

Let $X$ be a smooth projective variety. Consider a (reduced) simple normal corssing (SNC) divisor of $X$: $$D=\sum_{i=1}^sD_i.$$

\subsection{Log bend and break}
The main result of this paper is the following:

\begin{theorem}\label{thm:main}
Let $(X,D)$ be as above. Let $C\subset X$ be an integral complete curve not contained in $D$. Suppose that
$$(K_X+D)\cdot C<0.$$
Then there is an open and dense subset $C^\circ\subset C\setminus D$ such that through any point of $C^\circ$ there is a rational curve $R\subset X$, not contained in $D$, such that its normalization $\nu:\PP^1\to R\subset X$ satisfies
\begin{equation}\label{eq:lbb-contact}
 \#\Supp(\nu^*D)\leq1.
\end{equation}
In other words, $R$ intersects the divisor $D$ in at most $1$ point.
\end{theorem}

Theorem \ref{thm:main} proves\footnote{As already noted in \cite{KM}, Theorem \ref{thm:main} provides a simpler approach to some results of \cite{KM}.} the conjecture of Keel--McKernan \cite[Definition~1.2 and Conjecture~1.11]{KM}. The conjecture is termed {\em log bend and break}, because the case $D=\emptyset$ was proven by the deformation theory and mod-$p$ techniques collectively known as {\em bend and break}. 

Proposition \ref{prop:ample-bound} gives an explicit bound on the degree of $R$ in Theorem \ref{thm:main} depending on the chosen projective embedding of $X, D, C$ used in the proof.

Over the years, we have tried to prove Theorem \ref{thm:main} by extending the bend-and-break argument: when $D$ is a smooth divisor, we tried to use J. Li's space of stable relative maps \cite{Li}; for general $D$, we tried to use the spaces of stable log maps with maximal contacts \cite{C}, \cite{AC}, \cite{GS}. Unfortunately, we were not successful: while bend-and-break argument did produce rational curves in the situation of Theorem \ref{thm:main}, we were not able to make condition (\ref{eq:lbb-contact}) hold for such rational curves. With the recent advance in a new proof of Miyaoka-Mori criterion \cite{GX}, and with the help of GPT-6 Astra, we adopt the constructions of \cite{GX} to the case of log pair $(X,D)$ and derive Theorem\ref{thm:main}. We give an outline of the construction below.

\noindent{\it From logarithmic degree to an ordinary smooth family}

Let $f:B\to C\subset X$ be the normalization of the given curve. Choose an embedding 
\begin{equation}\label{eqn:embedding}
X\hookrightarrow\PP^N
\end{equation}
 and put
\begin{equation}\label{eqn:some_numbers}
 A=f^*\cO_X(1),\quad a=\deg A>0,\quad \delta=-(K_X+D)\cdot C>0,\quad r=\dim X+1.
\end{equation}
The method of \cite{GX} starts with weighted jets and obtains polynomial maps by controlling their coefficients. We use that method for a family which records information about the divisor $D$ at every point of $B$.

Lift $f$ to the nonzero twisted affine cone $\mathcal Z^\times$ of the embedding (\ref{eqn:embedding}), with maps $ (p,q):\mathcal Z^\times\longrightarrow X\times B$. For each component $D_i$ of $D$, consider the universal isomorphism 
$u_i:p^*\mathcal{O}_X(D_i)\xrightarrow{\ \sim\ }q^*(f^*\mathcal{O}_X(D_i))$ on the Isom scheme $\text{Isom}(p^*\mathcal{O}_X(D_i), q^*(f^*\mathcal{O}_X(D_i)))$ and impose
$$u_i(p^*d_i)=q^*(f^*d_i),$$
where $d_i$ is the canonical section of $\cO_X(D_i)$. These equations define, in the fiber product of Isom schemes $\text{Isom}(p^*\mathcal{O}_X(D_i), q^*(f^*\mathcal{O}_X(D_i)))$, a smooth family $W\to B$ near a distinguished section $\widehat s$. In local coordinates, they are $e_i z_i=b_i$, with $e_i$ invertible. 
The relative tangent bundle along the section is a vector bundle $E$ with
\begin{equation}\label{eqn:more_numbers}
 \rk E=r,\qquad \det E\simeq f^*\cO_X(-K_X-D), \qquad \deg E=\delta.
\end{equation}
In Section~\ref{sec:frames}, we give this construction and the calculation of the determinant. 

\noindent{\it From weighted jets to a ruled surface}

In Section~\ref{sec:jets}, we consider the scheme of order-$k$ jets to $W\to B$ and prove formulas about intersection numbers on its weighted projectivization (wher we assign weighted $q$ to the order-$q$ coefficient of a based jet):
\[
 H^{rk-1}F=(k!)^{-r},\qquad H^{rk}=-(k!)^{-r}\delta h_k, \qquad h_k=\sum_{q=1}^k\frac1q.
\]
The negative top intersection gives a curve of negative weighted degree. After a finite cover $\rho:B'\to B$, its weighted coefficients give a based jet with parameter line bundle $L$ satisfying
\[
 \mu:=\frac{\deg L}{\deg\rho}>\frac{\delta h_k}{2rk}>0.
\]
The sign is essential: $L$ is the inverse of the weighted tautological line along the chosen curve.

The cone coordinates of this jet have coefficients in $H^0(B',\rho^*A\otimes L^{-q})$. Hence the coefficient of order $q$ vanishes whenever $q>a/\mu$. Choose $\Delta>0$ that bounds the degrees of
homogeneous generators of the ideals of $X\subset\mathbb{P}^N$ and every $D_i\subset\mathbb{P}^N$, and take $k$ large enough so that the following holds:
\[
 \delta h_k>2r\Delta a.
\]
Every substituted defining equation then has degree less than $k$. An identity modulo $t^{k+1}$ is consequently an exact identity. This gives a rational map
\[
 \Phi:\Tot_{B'}(L)\dashrightarrow X, \qquad \Phi|_\sigma=f\rho,  
\]
regular near the entire zero section $\sigma$, where $\rho: B' \to B$ is a base change of curves.

More precisely, with $T_i=(f\rho)^*D_i$, we prove in Section~\ref{sec:polynomial} that 
\begin{equation}\label{eqn:divisor_equality}
 \Phi^*D_i=\pi_L^*T_i
\end{equation}
as Cartier divisors on a neighborhood of that entire section, see Lemma \ref{lem:boundary}. This is a new improvement relative to \cite{GX}, as this plays a key role in controlling how rational curves constructed in the ruled-surface argument meet the divisor components. 
The same positivity of $L$ excludes a constant generic projective fiber: scalar coefficients would be sections of $L^{-q}$ and would all vanish.

We then compactify $\Tot(L)$ to the space $\PP_{\mathrm{lines}}(\cO_{B'}\oplus L)$ of lines, with infinity section $\tau$: this is a ruled surface. The identity
\[
 [\sigma]=[\tau]+\pi_L^*c_1(L)
\]
shows that the only horizontal integral curve disjoint from $\sigma$ is $\tau$. We resolve $\Phi$ away from its neighborhood of $\sigma$. All resulting exceptional curves are vertical in $\PP_{\mathrm{lines}}(\cO_{B'}\oplus L)\to B'$. By (\ref{eqn:divisor_equality}), a pull-back of $D_i$ that is horizontal in $\PP_{\mathrm{lines}}(\cO_{B'}\oplus L)\to B'$ cannot meet $\sigma$, so every horizontal pull-backs of $D_i$ is at infinity. A general ruling fiber gives the required rational curve. This argument is presented in details in Section~\ref{sec:ruled}.

\subsection{Uniruledness of open locus}
We present an application of Theorem \ref{thm:main} to rational curves on (log) Fano surfaces.

\begin{proposition}[c.f. Theorem 1.3 and 1.3.1 of \cite{KM}]\label{prop:logbb-quotient-surfaces}
Let $S$ be a normal projective surface with quotient singularities, and let $D$ be a reduced effective Weil divisor on $S$, possibly empty. Suppose that
\begin{equation}\label{eq:logbb-negative-polarization}
(K_S+D)\cdot H<0
\end{equation}
for some ample Cartier divisor $H$ on $S$. Then the following assertions hold for every finite set
$V\subset S$.
\begin{enumerate}
\item
There is a dominating family of complete rational curves $R\subset S$, not contained in $D$, such that their normalization maps $\nu_R:\mathbb P^1\to S$ satisfy
\begin{equation}\label{eq:logbb-one-contact}
    \#\nu_R^{-1}\bigl(D\cup\operatorname{Sing}(S)\cup V\bigr)\leq 1.
\end{equation}

\item
If $(S,D)$ is purely log terminal at every point of $\operatorname{Sing}(S)$, then there is a dominating family whose members satisfy
\begin{equation}\label{eq:logbb-smooth-locus}
    R\subset S^{\mathrm{sm}},
    \qquad
    R\cap(V\setminus D)=\varnothing,
    \qquad
    \#\nu_R^{-1}(D)\leq 1.
\end{equation}
In particular, these curves meet $D\cup V$ at most once.
\end{enumerate}
No assumption on the Picard number of $S$ is required.
\end{proposition}
The proof of this proposition is given in Section \ref{sec:uniruled}.

\begin{corollary}\label{cor:log-del-pezzo-with-boundary}
Let $S$ be a normal projective surface with quotient singularities, and let $D$ be a reduced divisor. Suppose that $-(K_S+D)$ is ample. For every finite set $V\subset S$, $S$ is dominated by rational curves meeting $D\cup V$ at most once. If, in addition, $(S,D)$ is purely log terminal
along $\operatorname{Sing}(S)$, these rational curves may be chosen entirely in $S^{\mathrm{sm}}$ and disjoint from $V\setminus D$.

In particular, if $D=\emptyset$ and $-K_S$ is ample, then $S^{\mathrm{sm}}$ is dominated by complete rational curves avoiding any prescribed finite subset.
\end{corollary}

\begin{proof}
For any ample Cartier divisor $H$, the ampleness of $-(K_S+D)$ gives
\[
(K_S+D)\cdot H<0.
\]
The Corollary follows from Proposition~\ref{prop:logbb-quotient-surfaces}. When $D=\emptyset$, the required local purely log terminal condition holds because quotient surface singularities are klt.
\end{proof}

\begin{remark}\label{rem:comparison-keel-mckernan}
The first statement in Corollary \ref{cor:log-del-pezzo-with-boundary} extends the formulation of \cite[Theorem~1.3]{KM}: neither Picard number one nor a log-terminal condition on the reduced-divisor pair is needed for the stated divisor-contact conclusion. Under the original log-terminal hypothesis, the pair is purely log terminal at singular points of the surface;
see \cite[Appendix~L]{KM}. The second statement\footnote{We have also tried unsuccessfully to prove the second statement using bend-and-break technique on orbifold surfaces.} therefore also extends the smooth-locus formulation
in \cite[Theorem~1.3.1]{KM}.

\end{remark}

\begin{remark}[On plt assumption]\label{rem:boundary-condition-essential}
The additional hypothesis in the smooth-locus statement cannot simply be omitted, even when 
$-(K_S+D)$ is ample and the pair is log canonical.

Here is an example. Let
\[
    S=\mathbb P(1,1,m),\qquad m\geq 2,
\]
with weighted homogeneous coordinates $[x:y:z]$, and put
\[
    D=D_x+D_y,\qquad  D_x=\{x=0\},\qquad D_y=\{y=0\}.
\]
Writing $L$ for the ample $\mathbb Q$-Cartier class $\mathcal O_S(1)$, we have
\[
    K_S\sim_{\mathbb Q}-(m+2)L,\qquad
    D\sim_{\mathbb Q}2L,\qquad
    -(K_S+D)\sim_{\mathbb Q}mL.
\]
The divisors $D_x$ and $D_y$ are ample
$\mathbb Q$-Cartier divisors, and
\[
    D_x\cap D_y=\{v\},
    \qquad v=[0:0:1]\in\operatorname{Sing}(S).
\]
Every complete curve in $S^{\mathrm{sm}}$ must meet both $D_x$ and $D_y$. Since it avoids $v$, these
intersections occur at distinct points. Such a curve therefore meets $D$ at least twice on its normalization.

Locally at $v$, the pair is the quotient of the two coordinate axes in $\mathbb A^2$ by the scalar action of $\boldsymbol{\mu}_m$. It is log canonical but not purely log terminal there. Thus the divisor-contact conclusion on $S$ remains valid, while its proposed strengthening inside $S^{\mathrm{sm}}$ fails.
\end{remark}

\subsection{Campana curves}
We also consider an application of Theorem \ref{thm:main} to Campana's program \cite{C0}, \cite{C1}. Let $(X,D)$ be as above. 
Write $D=\sum_{i=1}^sD_i$ and, for each $i$, choose a positive integer $m_i$. Put 
\begin{equation}\label{cp:eq:pair}
\De=\sum_{i=1}^s\left(1-\frac1{m_i}\right)D_i,
\end{equation}
and write
\begin{equation}\label{cp:eq:L}
 \mathcal H=-(K_X+\De).
\end{equation}
We call the klt pair $(X,\De)$ {\em a smooth SNC Campana orbifold}, and {\em a weak Fano orbifold} when $\mathcal H$ is big and nef.

\begin{definition}\label{cp:def:campana}
A {\em Campana rational curve} on $(X,\De)$ is a nonconstant morphism $g:\PP^1\to X$ whose image is not contained in $D$ and such that every coefficient $c_{a,i}$ of $a\in \mathbb{P}^1$ in $g^*D_i$ satisfies
\[
 c_{a,i}=0\quad\hbox{or}\quad c_{a,i}\geq m_i.
\]
It is {\em divisible} if $m_i\ |\ c_{a,i}$ for every $a,i$. 

The orbifold $(X,\De)$ is {\em Campana uniruled} if a finite-type algebraic family of Campana rational curves has dominant evaluation. It is {\em Campana rationally connected} if such a family has dominant two-point evaluation to $X\times X$.
\end{definition}

Multiple covers are allowed in this definition. Thus divisibility is a property of the map, and need not hold for the normalization map of its image. 

Our application of Theorem \ref{thm:main} is the following
\begin{theorem}[Free Campana curves]\label{cp:thm:main}
For the log pair $(X,\De)$, suppose that $\mathcal H$ is big. Then there is a covering family of free\footnote{Here, freeness means logarithmic freeness with the contact points marked, in the sense of Chen--Lehmann--Tanimoto \cite{CLT}. Its cohomological definition and its relation to the normal bundle of the graph are given in Section~\ref{cp:sec:free}.} divisible Campana rational curves with fixed numerical class and fixed contact data. 
\end{theorem}
Theorem \ref{cp:thm:main} and \cite[Proposition 4.4]{CLT} imply that $(X,\De)$ is Campana uniruled. Campana conjectured  that a smooth projective orbifold pair $(X, \De)$ over $\mathbb{C}$ is orbifold uniruled if and only if $K_X+\De$ is not pseudo-effective \cite[Conjecture~11.5(1)]{Ca}. Our result can be viewed as a piece of evidence for Campana's conjecture.

\begin{corollary}\label{cp:cor:rhoone}
In the situation of Theorem~\ref{cp:thm:main}, if $\rho(X)=1$, then $(X,\De)$ is Campana rationally connected.
\end{corollary}

In the future work \cite{CT}, we will study Campana rational connectedness conjecture \cite[Conjecture 1.2]{LT} in the case of higher Picard rank.

Theorem~\ref{cp:thm:main} does not require nefness. In the Picard number one case, bigness of $\mathcal H$ implies ampleness, and rational connectedness of $X$ follows. The connecting curves given by Corollary~\ref{cp:cor:rhoone} satisfy the Campana inequalities.

We now give an outline of the proof of Theorem \ref{cp:thm:main}. Given $(X,\De)$, we construct in Section~\ref{cp:sec:aux} a normal projective variety $Y$, a projective map $f:Y\to X$, and a reduced boundary divisor $B_Y\subset Y$ such that
\[
 K_Y+B_Y\sim_{\Q}f^*(K_X+\De),\qquad  U=Y\setminus B_Y\text{ is smooth},\qquad f^*D_i=m_iE_i.
\]
The fibers of $U\to X$ over $X\setminus D$ are tori. The construction uses separate cyclic covers of a product of $\PP^1$-bundles and keeps both section divisors. It is adapted from the auxiliary construction of Lehmann--Tanimoto \cite[Section~3]{LT}.

A general ample complete intersection curve on $Y$ avoids its codimension-$2$ non-SNC locus. Its strict transform on a log resolution has negative logarithmic canonical degree, because $\mathcal H$ is big. Theorem~\ref{thm:main} gives rational curves with at most one boundary point. Their projections to $X$ are nonconstant: a vertical curve would give a nonconstant map $\A^1\to(\Gm)^s$. All projected contacts are divisible except possibly at that one point. A self-map of $\PP^1$ of degree $\lcm(m_1,\ldots,m_s)$, totally ramified there, modifies the remaining contacts as needed.

To make this pointwise argument into a statement about families, we use finite-type incidence schemes and the uncountability of $\C$. 
After fixing the contact data, the derivative of the evaluation map produces global sections of the normal bundle of the graph; its general member is therefore free. Finally, in Picard number one case, the common nonzero covering class is interior to the nef cone of curves. The rational connectedness result of Chen--Lehmann--Tanimoto \cite[Corollary~6.7]{CLT} then applies.

The rest of this paper is divided into two Parts. Part 1, which consists of Sections~\ref{sec:jets}--\ref{sec:ruled}, contains a proof of Theorem \ref{thm:main}. Proposition \ref{prop:logbb-quotient-surfaces} is proved in Section \ref{sec:uniruled}. Part 2, which consists of Sections~\ref{cp:sec:aux}--\ref{cp:sec:free}, establishes consequences on Campana curves. 

\subsection{Acknowledgment}
In the process of developing the results presented here and writing this article, the authors greatly benefited from the assistance of GPT-6 Astra. J.-C. C. is supported by the National Science and Technology Council (NSTC) of Taiwan under Grant No. NSTC-114-2115-M-007-004- and NSTC-115-2115-M-007-013-.

{\large \part{Existence of rational curves}}

In this part, we present a proof of Theorem \ref{thm:main}.

\section{Weighted jets of an ordinary smooth family}\label{sec:jets}

Let $B$ be a smooth projective connected curve. Consider a separated morphism $q$ of finite type, with a section $s$: 
$$q:W\to B,\quad s:B\to W.$$
Assume that $q$ is smooth of relative dimension $r\geq1$ on a
neighborhood of $s(B)$. Hence we may replace $W$ by this neighborhood and consider:
\[
 E=s^*T_{W/B},\qquad \deg E=\delta.
\]
The purpose of this section is to construct jets of $W$ whose parameters are in a line bundle of positive degree when $\delta>0$. 

\subsection{The based jet space and its coordinates}

Fix an integer $k\geq2$. For a $B$-scheme $T$, put
\[
 T[k]=\Spec_T\bigl(\cO_T[t]/(t^{k+1})\bigr).
\]
The based $k$-jet functor assigns to $T$ the set of $k$-jets
\[
 J_k(W,s)(T)=
 \bigl\{j:T[k]\longrightarrow W\text{ over }B
 \, \big| \, j|_{t=0}=s|_T\bigr\}.
\]
We first describe the scheme representing this functor. We refer to \cite{EM} for basics about jet schemes.

Fix $b\in B$. After restricting to an open neighborhood $U\subset B$ of
$b$, there is an open neighborhood $V\subset W_U$ of $s(U)$ and an
\'etale $U$-morphism
\[
 (z_1,\ldots,z_r):V\longrightarrow\A^r_U.
\]
Subtracting the values $s^*z_i$, we may assume that $z_i\circ s=0$.
For $T\to U$, a based map $T[k]\to\A^r_U$ is specified uniquely by
\[
 z_i(t)=\sum_{q=1}^k x_{iq}t^q,
 \qquad x_{iq}\in H^0(T,\cO_T),\quad 1\leq i\leq r.
\]
Formal \'etaleness gives a unique lift to $V$ with prescribed reduction
$s|_T$. Conversely, every based jet into $W_U$ factors through $V$:
the underlying topological space of $T[k]$ is that of $T$, and its image
is contained in $s(U)\subset V$. Thus
\begin{equation}\label{eqn:chart_jet_scheme}
 J_k(W,s)|_U\simeq\A_U^{rk}.
\end{equation}
The uniqueness of the lifts makes these identifications compatible on
overlaps; they therefore represent the functor by an affine $B$-scheme.
This uses only \'etale coordinates near the section, not a product
decomposition of $W$; see also~\cite[Tags 054L and 02HF]{Stacks}.

Write $x_q=(x_{1q},\ldots,x_{rq})$. For $\lambda\in \mathbb{G}_m$, precomposing with $t\mapsto\lambda t$ defines a $\mathbb G_m$-action on the jet scheme for which
\begin{equation}\label{eqn:Gm-action}
 \lambda\cdot x_{iq}=\lambda^q x_{iq}.
\end{equation}
Therefore, $x_{iq}$ are assigned {\em weight} $q$. 

We describe a change of coordinates using the formal completions along
the section. \'Etale coordinates identify these completions with formal
affine space over $U$. A second centered coordinate system is therefore
a formal power series in the first, with zero constant term and
invertible linear term. Substituting $\sum_{q=1}^k x_qt^q$ and taking
the coefficient of $t^q$ gives
\begin{equation}\label{eq:triangular}
 x'_q=g x_q+P_q(x_1,\ldots,x_{q-1}),
\end{equation}
where $g$ is the change of components for vectors in $E$, and $P_q$ is
weighted homogeneous of weight $q$. Every monomial of $P_q$ has ordinary
degree at least two. Indeed, a monomial containing a variable of weight
$q$ and at least one other variable would have weight greater than $q$.

\subsection{Weighted projectivization}

Consider the sheaf of polynomial functions on $J_k(W,s)$, with grading given by weights:
$$\mathcal A=\bigoplus_{\ell\geq0}\mathcal A_\ell.$$ 
In the chart (\ref{eqn:chart_jet_scheme}), we have
\[
 \mathcal A|_U=\cO_U[x_{iq}:1\leq i\leq r,\ 1\leq q\leq k],
 \qquad \operatorname{wt}(x_{iq})=q.
\]
Weighted homogeneity in \eqref{eq:triangular} makes this grading independent of the chart. Each $\mathcal A_\ell$ is locally free: the monomials of weight $\ell$ form a finite local basis. Define
\[
 Y_k=\Proj_B\mathcal A\xrightarrow{\,\pi\,} B,\quad m=rk.
\]
Locally over $B$, the morphism $\pi$ is the product with the weighted projective space having $r$ coordinates of each weight $1,\ldots,k$. Fibers of $\pi$ have dimension $rk-1$, so $$\dim Y_k=m.$$ Alternatively, $Y_k$ is the weighted projective quotient of the scheme $J_k(W,s)$ of nonzero based jets by the $\mathbb{G}_m$-action in (\ref{eqn:Gm-action}). In particular, over an algebraically closed field, a point is represented by a nonzero tuple $(x_{iq})$, up to $x_{iq}\mapsto\lambda^q x_{iq}$.

Choose $M$ divisible by $1,\ldots,k$ and sufficiently large, so that the Veronese algebra $\mathcal A^{(M)}$ is generated in degree $1$. Then the surjection
\[
 \Sym\mathcal A_M\twoheadrightarrow\mathcal A^{(M)}
\]
defines a closed immersion $$Y_k\hookrightarrow\PP_{\mathrm{quot}}(\mathcal A_M)$$
over $B$. Here $\PP_{\mathrm{quot}}(V)=\Proj_B\Sym V$ parametrizes line bundle quotients of $V$; its tautological bundle is the universal quotient $\cO(1)$. Thus $$\cO_{Y_k}(M)$$ is an actual line bundle and is relatively very ample; see~\cite[Tag 01VL]{Stacks}. Since $B$ is projective, $Y_k$
is also projective. It can be checked on the usual weighted polynomial Proj charts that $Y_k$ is normal. Consider the $\mathbb{Q}$-divisor class
\begin{equation}\label{eqn:divisor-H}
 H=\frac1M c_1\bigl(\cO_{Y_k}(M)\bigr)\in \text{Div}(Y_k)_\mathbb{Q},
\end{equation}
and let $$F\in \text{Div}(Y_k)$$ 
be the numerical divisor class of a fiber of $\pi: Y_k\to B$. Passing to a further divisible Veronese replaces $\cO(M)$ by the corresponding tensor power, so $H$ is independent of this choice.

\begin{remark}[Sign convention] 
In weight one,
$Y_1=\Proj_B\Sym E^*$ parametrizes lines in $E$. A line $L\subset E$ determines the quotient $E^*\twoheadrightarrow L^{-1}$. Consequently the pullback of $\cO_{Y_1}(1)$ is $L^{-1}$, not $L$.
The same sign will occur in every weight.
\end{remark}

The following lemma is about two intersection numbers on $Y_k$ that play an important role in the construction.

\begin{lemma}[Intersection numbers]\label{lem:intersection}
With $h_k=\sum_{q=1}^k 1/q$ the $k$-th harmonic number, $v=(k!)^{-r}$, and $\delta=\text{deg}\,E$, we have
\begin{equation}\label{eq:intersection}
 H^{m-1}\cdot F=v,\qquad H^m=-v\delta h_k.
\end{equation}
\end{lemma}

\begin{proof}
(\ref{eq:intersection}) is the analogue of the results in \cite[Proposition 2.4]{GX}. We compute ranks and degrees by a finite filtration, so that no deformation of the family $W$ is needed.

\noindent\emph{The filtration.}
Let $I\subset\mathcal A$ be the ideal of the constant jet, and put
\[
 F^p\mathcal A_\ell=(I^p)_\ell.
\]
In a jet chart, $F^p\mathcal A_\ell$ is spanned by the monomials of weight $\ell$ and ordinary degree at least $p$. The coordinate changes \eqref{eq:triangular}, as well as their inverses, preserve this submodule. Since every variable has weight at least $1$, $F^{\ell+1}\mathcal A_\ell=0$. On the associated graded module, all the nonlinear terms disappear, leaving the coordinate changes of $E^*$ in each weight block. Here the coefficients $x_q$ transform as components of tangent vectors, whereas linear polynomial functions of these coefficients belong to the dual bundle. In particular, the weight-$q$ summand of $I/I^2$ is $E^*$. Hence
\[
 \operatorname{gr}^p_F\mathcal A_\ell\simeq
 \bigoplus_{(a_q)\in\mathsf{I}_{p,\ell}} \bigotimes_{q=1}^k\Sym^{a_q}E^*, 
\quad \text{where } \mathsf{I}_{p,\ell}=\Big\{(a_q)\, \Big|\,\sum_{q=1}^k a_q=p, \sum_{q=1}^k q a_q=\ell \Big\}.
\]
The monomial description shows that $F^p\mathcal A_\ell$ and its successive quotients are locally free. Thus $F^\bullet\mathcal{A}_\ell$ is a finite filtration by vector subbundles, and both rank and degree are additive along it. In particular, they agree with the rank and degree of the weight-$\ell$ part of
\[
 \Sym\left(\bigoplus_{q=1}^k E^*_{(q)}\right),
\]
where $E^*_{(q)}\subset E^*$ is the weight-$q$ summand.
 
\noindent\emph{Ranks and degrees of the graded pieces.}
For a rank-$r$ bundle $V$, let $c_a=\binom{a+r-1}{r-1}=\rk(\Sym^a V)$. The following formula can be proved by splitting principle:
\[
 \deg\Sym^a V=\frac{a}{r}c_a\deg V.
\]

The generating series for ranks of graded pieces of $\Sym^\bullet(E_{(q)}^*)$ is
\[
 R_q(z):=\sum_{a\geq 0} (\text{rk}\Sym^a E_{(q)}^*)z^{\text{wt}\Sym^a E_{(q)}^*}=\sum_{a\geq0}c_a z^{qa}=(1-z^q)^{-r}.
\]
Since $\deg E^*=-\delta$, the generating series for degrees of graded pieces of $\Sym^\bullet(E_{(q)}^*)$ is
\[
 \begin{aligned}
 D_q(z)&:=\sum_{a\geq 0} (\text{deg}\Sym^a E_{(q)}^*)z^{\text{wt}\Sym^a E_{(q)}^*}\\
  &= -\frac{\delta}{r}\sum_{a\geq0}a c_a z^{qa}\\
  &= -\delta\frac{z^q}{(1-z^q)^{r+1}}
   =-\delta R_q(z)\frac{z^q}{1-z^q}.
 \end{aligned}
\]
The middle equality follows by applying $u\,d/du$ to $(1-u)^{-r}$ and then putting $u=z^q$. For tensor products $U\otimes V$, we have $\text{rk}(U\otimes V)=\text{rk}(U)\cdot \text{rk}(V)$ and
\[
 \deg(U\otimes V)=\rk(V)\deg U+\rk(U)\deg V.
\]
Consequently, the generating series for degrees of graded pieces of $\mathcal{A}$ is obtained by selecting one degree factor $D_q$ and rank factors in all other blocks. Writing
$P(z)=\prod_{q=1}^k(1-z^q)^{-r}$, we obtain
\begin{equation}\label{eq:ranks}
 \sum_{\ell\geq0}\rk(\mathcal A_\ell)z^\ell =P(z),
\end{equation}
\begin{equation}\label{eq:degrees}
 \sum_{\ell\geq0}\deg(\mathcal A_\ell)z^\ell
=\sum_{q=1}^kD_q(z)\prod_{j\ne q}R_j(z)
=-\delta P(z)\sum_{q=1}^k\frac{z^q}{1-z^q}.
\end{equation}

\smallskip\noindent
\emph{Leading coefficients.}
Since $1-z^q=q(1-z)+O((1-z)^2)$, the expansions at $z=1$ are
\[
 \begin{aligned}
 P(z)&=v(1-z)^{-m}+O\bigl((1-z)^{-m+1}\bigr),\\
 \sum_{q=1}^k\frac{z^q}{1-z^q}
   &=h_k(1-z)^{-1}+O(1).
 \end{aligned}
\]
All remaining poles of these as functions of $z$ are roots of unity. If $\zeta\ne1$ has order $d$, put $$\nu=\#\{q\, |\, 1\leq q\leq k,\ d\mid q\}.$$ 
The pole order of $P(z)$ at $z=\zeta$ is $r\nu\leq m-1$, because the weight-one block is regular there. The pole order of the function in (\ref{eq:degrees}) is at most $r\nu+1\leq m$. The coefficient of $z^\ell$ in $(1-z)^{-a}$ is $\binom{\ell+a-1}{a-1}$; a pole of order $a$ at another root of unity contributes a bounded root-of-unity factor times a polynomial of degree
at most $a-1$. Taking partial fractions then gives
\[
 \begin{aligned}
 \rk\mathcal A_\ell
   &=\frac{v}{(m-1)!}\ell^{m-1}+O(\ell^{m-2}),\\
 \deg\mathcal A_\ell
   &=-\frac{v\delta h_k}{m!}\ell^m+O(\ell^{m-1}).
 \end{aligned}
\]
Here $m\geq2$, since $k\geq2$.

\noindent\emph{Passage to intersection numbers.}
Take $\ell=Ma$ with $a$ sufficiently large. By the standard Proj identities and relative Serre vanishing, we have
\[
 \pi_*\cO_{Y_k}(\ell)=\mathcal A_\ell,\qquad  R^j\pi_*\cO_{Y_k}(\ell)=0\quad j>0.
\]
Write $g_B$ for the genus of $B$, then by Riemann--Roch, we obtain
\[
 \begin{aligned}
 \chi(Y_k,\cO_{Y_k}(\ell))
  &=\chi(B,\mathcal A_\ell)\\
  &=\deg\mathcal A_\ell+(1-g_B)\rk\mathcal A_\ell\\
  &=-\frac{v\delta h_k}{m!}\ell^m+O(\ell^{m-1}).
 \end{aligned}
\]
On the other hand, because $Y_k$ is projective of dimension $m$, the Euler polynomial $\chi(Y_k,\cO_{Y_k}(\ell))$ of the line bundle $\cO_{Y_k}(M)$ as $a^m$-coefficient equal to $$\frac{c_1(\cO_{Y_k}(M))^m}{m!}=\frac{(MH)^m}{m!}.$$ 
The case $\ell=Ma$ then implies $H^m=-v\delta h_k$.

On a fiber $Y_{k,b}$ over $b\in B$, the dimension of the weight-$\ell$ polynomial piece is $\rk\mathcal A_\ell$. For large and divisible $\ell$, this is equal to the Euler characteristic of $\cO_{Y_{k,b}}(\ell)$, by the same Proj calculation. Comparing its leading coefficient, now in dimension $m-1$, we obtain $H^{m-1}\cdot F=v$.
\end{proof}

\subsection{Construction of a positive parameter line bundle}
The following lemma constructs a suitable family of based jets. 

\begin{lemma}[A positive parameter line bundle]\label{lem:positive}
Assume $\delta>0$. For every $k\geq2$, after a finite surjective map of smooth projective connected curves $\rho:B'\to B$, there are a line bundle $L$ on $B'$ and a based $k$-jet
\begin{equation}\label{eq:Ljet}
 j:\Spec_{B'}\left(\bigoplus_{q=0}^k L^{-q}\right)\longrightarrow W\times_B B'
\end{equation}
whose coefficients do not all vanish at any point, and for which
\begin{equation}\label{eq:slope}
 \mu:=\frac{\deg L}{\deg\rho}>\frac{\delta h_k}{2rk}>0.
\end{equation}
The algebra in \eqref{eq:Ljet} has the truncated multiplication; its source is the $k$th infinitesimal neighborhood of the zero section in $\Tot(L)$.
\end{lemma}

\begin{proof}
\emph{A curve of negative weighted degree.}
Set $$\theta=\delta h_k/(2m).$$ Since a fiber $F$ of $\pi: Y_k\to B$ satisfies $F^2=0$, Lemma~\ref{lem:intersection} gives
\[
 (H+\theta F)^m =H^m+m\theta H^{m-1}\cdot F =-\tfrac12v\delta h_k<0.
\]
A nef $\mathbb{Q}$-Cartier divisor on a projective variety has nonnegative top self-intersection. Thus $H+\theta F$ is not nef, and there is an integral curve 
$$\Gamma\subset Y_k, \quad \text{with } (H+\theta F)\cdot\Gamma<0.$$ 
A curve contained in a fiber of $\pi$ has positive $H$-degree, because $H$ is $\pi$-relatively ample. On the other hand, the intersection of such a curve with a fiber $F$ is $0$. Therefore, $\Gamma$ dominates $B$. Let $d_\Gamma$ be the degree of the map $\Gamma\to B$, then $F\cdot\Gamma=d_\Gamma$, and hence
\begin{equation}\label{eq:gamma}
\frac{-H\cdot\Gamma}{d_\Gamma}>\theta.
\end{equation}

\noindent\emph{An affine representative over a finite extension.}
We now construct the parameter line bundle, corresponding to
\cite[Lemma 3.1 and Proposition 3.2]{GX}. Let $K=\C(\Gamma)$. The generic point of $\Gamma$ is a $K$-point of $Y_k$ above the generic point of $B$. The nonzero affine weighted cone over this point has a point over an algebraic closure of $K$, hence is a nonempty scheme of finite type over $K$. A closed point of this scheme has residue field being a finite extension $K_1/K$. It gives an affine representative of the generic weighted jet over $K_1$.

Let $C_1$ be the normalization of $\Gamma$ in $K_1$. This is a smooth projective connected curve, finite over $\Gamma$ and hence over $B$. In a based-jet chart, the chosen representative has coefficients $x_{iq}\in K_1$, not all $0$. On another chart, these coefficients are related by \eqref{eq:triangular}. They need not yet be regular or nowhere simultaneously zero along $C_1$.

\noindent\emph{Weighted valuations.}
For $y\in C_1$, choose a jet chart above its image in $B$ and put
\[
 \alpha_y=\min_{i,q}\frac{\ord_y(x_{iq})}{q}\in \mathbb{Q},
\]
omitting the zero coefficients. This rational number is independent of the chart. Indeed, every monomial of weight $q$ in \eqref{eq:triangular} has valuation at least $q\alpha_y$, since its
coefficients are regular at $y$. Thus the new coefficients satisfy $\ord_y(x'_{iq})\geq q\alpha_y$, and the minimum computed in the new chart is at least $\alpha_y$. Applying the same argument to the inverse coordinate change gives the opposite inequality.

Only finitely many $\alpha_y$ are nonzero. To see this, use one jet chart over a dense open subset of $B$. Away from finitely many zeros and poles, all its rational coefficients are regular and one chosen nonzero coefficient is a unit. The omitted points of $B$, and their inverse images in $C_1$, are finite as well.

We next make all $\alpha_y$ integral. Let $S$ be the set of all $y\in C_1$ where $\alpha_y\neq 0$, and let $N$ be the least common multiple of denominators of nonzero $\alpha_y's$, with the empty least common multiple defined to be one. Thus $N=1$ when all the $\alpha_y$ are integral. If $N>1$, then $S$ is nonempty; choose a point $p_0\in C_1\setminus S$ and put $\Sigma=\sum_{y\in S}y$. Choose an integer $a$ such that $a-2\#S>2g(C_1)-2$, then by degree consideration, we have
\[
 H^1\bigl(C_1,\cO_{C_1}(a p_0-2\Sigma)\bigr)=0.
\]
The exact sequence for restriction to $2\Sigma$ then gives a surjection
\[
 H^0\bigl(C_1,\cO_{C_1}(a p_0)\bigr) \longrightarrow \bigoplus_{y\in S}\cO_{C_1,y}/\mathfrak m_y^2.
\]
Here $\cO_{C_1}(a p_0)$ has its canonical trivialization near $S$. Choose a uniformizer $\varpi_y$ at each $y\in S$ and prescribe the image to be $\varpi_y\bmod\mathfrak m_y^2$. A preimage is a rational function $h$ on $C_1$, regular off $p_0$, for which $\ord_y(h)=1$ for every $y\in S$. 
Adjoin a root of $T^N-h$ to $K_1$. At each prescribed point, this polynomial is Eisenstein; the extension has degree $N$ and is totally ramified there. Normalize $C_1$ in this field extension. If $y'$ lies over $y$ with ramification index $e_{y'}$, then the pulled-back coefficients satisfy
\[
 \min_{i,q}\frac{\ord_{y'}(x_{iq})}{q} =e_{y'}\alpha_y.
\]
This is integral above $S$, and is zero away from $S$, even if extra ramification occurs there. If $N=1$, no additional cover is needed. Denote the resulting smooth projective connected curve by $B'$, and the composite map to $B$ by 
$$\rho: B'\to B.$$
From now on $\alpha_y$ denotes the integral minimum on $B'$.

\noindent\emph{Regular coefficients and gluing.}
Define
\[
 D_\alpha=\sum_{y\in B'}\alpha_y y,
 \qquad L=\cO_{B'}(D_\alpha).
\]
As seen above, this sum is finite; the divisor is allowed to have negative coefficients. The canonical rational section $1$ of $L$ has divisor $D_\alpha$. On an open set with a regular section $\varepsilon$ of $L$, write $1=\gamma\varepsilon$ and put
\begin{equation}\label{eqn:a_iq}
 a_{iq}=\gamma^{-q}x_{iq}.
\end{equation}
At every point in this open set, we have
\[
 \ord_y(a_{iq})=\ord_y(x_{iq})-q\alpha_y\geq0.
\]
Thus all the $a_{iq}$ are regular. The minimum defining $\alpha_y$ is attained, so at least one $a_{iq}$ has order $0$ and is a unit at $y$. This proves the nonvanishing at every point, as claimed.

The local coefficients define a based jet with coordinate expressions 
$z_i(t)=\sum_{q=1}^k a_{iq}t^q$. We check both changes needed for gluing. First, if $\varepsilon'=u\varepsilon$ is another regular section, then $\gamma'=u^{-1}\gamma$ and
\[
 a'_{iq}=u^q a_{iq},\qquad t'=u^{-1}t.
\]
Consequently, $a'_{iq}(t')^q=a_{iq}t^q$. Second, by (\ref{eqn:a_iq}), weighted homogeneity gives
\[
 \gamma^{-q}P_q(x_1,\ldots,x_{q-1}) =P_q(a_1,\ldots,a_{q-1}),
\]
so the normalized coefficients obey the same jet-coordinate transition rule \eqref{eq:triangular}. The local jets therefore agree on overlaps and define \eqref{eq:Ljet}. Its source is given in (\ref{eq:Ljet}) because the local parameter transforms as $t'=u^{-1}t$.

\noindent\emph{The tautological quotient and its degree.}
The nowhere simultaneously zero coefficients have a weighted projectivization, givin a map  $\tau:B'\to Y_k$. At the generic point, it is the original point of $\Gamma$, since multiplication by $\gamma^{-1}$ acts with weight $q$ on its $q$-th coefficients. Both maps are defined on all of $B'$, so separatedness of $Y_k$ shows that $\tau$ is the composite $B'\to\Gamma\hookrightarrow Y_k$.

To determine its tautological line bundle, we evaluate weighted polynomials of degree $M$:
\[
 \rho^*\mathcal A_M\longrightarrow L^{-M},\quad  R\mapsto R(a)\varepsilon^{-M}.
\]
This expression is independent of the regular section $\varepsilon$: replacing $\varepsilon$ by $u\varepsilon$ multiplies $R(a)$ by $u^M$ and $\varepsilon^{-M}$ by $u^{-M}$. It is also independent of the jet chart, since it is evaluation of a polynomial function on that jet. The map is surjective. At each point some $a_{iq}$ is a unit; as $q\,|\, M$, the monomial $x_{iq}^{M/q}$ evaluates to a unit in a local section of $L^{-M}$.

This line bundle quotient gives exactly the composite of $\tau$ with the Veronese embedding into $\PP_{\mathrm{quot}}(\mathcal A_M)$. Pulling back its universal quotient yields
\begin{equation}\label{eq:tautological}
 \tau^*\cO_{Y_k}(M)\simeq L^{-M},\quad  L^M\simeq\tau^*\cO_{Y_k}(-M).
\end{equation}
These are isomorphisms of line bundles, and fix the sign without a numerical convention. Let $e$ be the degree of $B'$ over the normalization of $\Gamma$. Taking degrees in \eqref{eq:tautological} and using $\tau_*[B']=e[\Gamma]$ gives
\[
 \deg L=-e(H\cdot\Gamma),\qquad  \deg\rho=e d_\Gamma.
\]
Dividing these identities and applying \eqref{eq:gamma}, we obtain
\[
 \frac{\deg L}{\deg\rho}
       =\frac{-H\cdot\Gamma}{d_\Gamma}
       >\theta=\frac{\delta h_k}{2rk}.
\]
This proves \eqref{eq:slope}, in particular $\deg L>0$.
\end{proof}

\section{An ordinary smooth family encoding the divisor}\label{sec:frames}

The purpose of this section is to construct a smooth family with a distinguished section such that the degree of its relative tangent bundle is $-(K_X+D)\cdot C$. We first lift the curve $C\subset X$ to an affine cone. We then introduce one invertible frame for each component of the boundary.
The frame equations subtract precisely the boundary degrees from the tangent degree.

\subsection{The twisted cone and its distinguished section}

Choose a very ample line bundle $\cO_X(1)$ and its associated embedding 
$$X\hookrightarrow\PP^N.$$
As before, denote by $f: B\to C\subset X$ the normalization of $C$. Put
\[
 A=f^*\cO_X(1),\qquad a=\deg A>0.
\]
The pullbacks via $f$ of the homogeneous coordinates of $\mathbb{P}^N$ are sections 
\begin{equation}\label{eqn:sections_of_A}
f_0,\ldots,f_N\in H^0(B,A)
\end{equation}
 with no simultaneous $0$.

Let 
\begin{equation}\label{eqn:twisted_cone}
\mathcal Z\subset\Tot_B(A^{\oplus(N+1)})
\end{equation}
be the twisted affine cone over $X$. Here is the meaning of this definition. If $F$ is a
homogeneous equation of $X$ of degree $e$, then substitution in the tautological tuple on $\Tot_B(A^{\oplus(N+1)})$ gives a section of the pullback of $A^e$. Its zero scheme is therefore well-defined, even though $A$ need not be trivial. The common zero scheme of all these sections is $\mathcal Z$. In a local regular section $\alpha$ of $A$, write a tuple as $(x_0\alpha,\ldots,x_N\alpha)$. If $\alpha'=v\alpha$, the scalar coordinates become $x'_j=v^{-1}x_j$,
and
\[
 F(x'_0,\ldots,x'_N)=v^{-e}F(x_0,\ldots,x_N).
\]
Thus the local cone equations agree on overlaps as line-bundle-valued equations.

Let $$\mathcal Z^\times=\mathcal{Z}\setminus (\text{zero section}).$$
Projectivization gives maps
\[
 p:\mathcal Z^\times\longrightarrow X, \quad q:\mathcal Z^\times\longrightarrow B.
\]
More intrinsically, a point of $\mathcal Z^\times$ over a scheme $T$ consists of maps $g:T\to X$, $b:T\to B$, and an isomorphism
\begin{equation}\label{eqn:line_bdle_isom}
 g^*\cO_X(1)\xrightarrow{\ \sim \ }b^*A.
\end{equation}
Indeed, such an isomorphism transports the homogeneous coordinate sections to a tuple in $b^*A$ generating that line bundle; conversely, a generating tuple determines the projective map and the isomorphism. It follows that
\[
 (p,q):\mathcal Z^\times\longrightarrow X\times B
\]
is the $\mathbb G_m$-torsor of isomorphisms between the two line bundles as in (\ref{eqn:line_bdle_isom}). The tuple $(f_0,\ldots,f_N)$ defines a section $s: B\to \mathcal{Z}^\times$ of $q$; in this description it is the pair $(f,\id_B)$ with the identity isomorphism $f^*\cO_X(1)=A$.

Since $X$ is smooth of dimension $n$, the morphism $q$ is smooth of relative dimension $n+1$. The relative tangent sequence of the $\mathbb{G}_m$-torsor, restricted to the section $s$, is
\begin{equation}\label{eq:cone-tangent}
 0\longrightarrow\cO_B\longrightarrow E_0:=s^*T_{\mathcal Z^\times/B}
 \longrightarrow f^*T_X\longrightarrow0.
\end{equation}
The first term is the Lie algebra of $\mathbb G_m$, which is a trivial line bundle. Taking determinants gives $\det E_0\simeq f^*\cO_X(-K_X)$, and hence
$$\deg E_0=-K_X\cdot C.$$

\subsection{Boundary frames and the smoothness calculation}

For each divisor component $D_i$, put
\[
 M_i=\cO_X(D_i),\qquad A_i=f^*M_i.
\]
Let $$d_i\in H^0(X,M_i)$$ 
be the canonical section with zero divisor $D_i$, and put 
$$b_i=f^*d_i\in H^0(B,A_i).$$ 
Since the curve $C$ is not contained in $D$, none of the $b_i$ is identically $0$.

Over $\mathcal Z^\times$, form the fiber product of frame torsors:
\[
 \mathcal T= \Isom(p^*M_1,q^*A_1)\times_{\mathcal{Z}^\times}\times_{\mathcal{Z}^\times}...\times_{\mathcal{Z}^\times} \Isom(p^*M_s,q^*A_s).
\]
Each factor parametrizes invertible isomorphisms of line bundles. Thus $\mathcal T\to\mathcal Z^\times$ is a $(\mathbb G_m)^s$-torsor. On $\mathcal T$, let 
$u_i:p^*M_i\xrightarrow{\ \sim\ }q^*A_i, i=1,...,s$ be the universal isomorphisms. Define $$W\subset \mathcal{T}$$ by the equations
\begin{equation}\label{eq:frames}
 u_i(p^*d_i)=q^*b_i,\qquad 1\leq i\leq s.
\end{equation}
The difference of the two sides is a section of $q^*A_i$; its zero scheme is therefore defined without choosing a regular section. As before, $p$ and $q$ also denote the composite maps from $\mathcal T$ and $W$. Along $s$, we have $s^*p^*M_i=A_i=s^*q^*A_i$. The identity isomorphisms give a lift $$\widehat s:B\to W.$$

\begin{lemma}\label{lem:W}
The morphism $W\to B$ is smooth of relative dimension $n+1$ in a neighborhood of $\widehat s$. The restriction of its relative tangent bundle to this neighborhood satisfies
\[
 \rk(\widehat s^*T_{W/B})=n+1,\qquad \deg(\widehat s^*T_{W/B})=\delta.
\]
Over $B\setminus f^{-1}D$, forgetting the frames identifies $W$ with the part of $\mathcal Z^\times$ lying over $X\setminus D$.
\end{lemma}

\begin{proof}
\noindent\emph{Smoothness.}
Fix $b\in B$, put $x=f(b)\in X$, and let $I=\{i\ |\ x\in D_i\}$. In a neighborhood of $x$, choose local sections of all $M_i$'s, and write $d_i=z_i$ in these sections. In a neighborhood of $b$, use their pullbacks as sections of $A_i$; then $b_i=z_i\circ f$. Write an isomorphism $u_i$ as multiplication by $e_i\in\mathbb G_m$. The identity isomorphisms have $e_i=1$, and the equations of $W$ become
\[
 e_i z_i-b_i=0.
\]
Their relative differentials at $\widehat s(b)$ are
\[
 dz_i+z_i(x)\,de_i.
\]
There is no derivative of $b_i$ here, since the differential is relative to $B$.

The simple normal crossings assumption says that the differentials $dz_i$, $i\in I$, are independent in $T_x^*X$. Choose tangent vectors $v_i\in T_xX$ dual to these differentials and lift them to the tangent space of $\mathcal Z^\times$ over $B$. For each $j\notin I$, use the tangent vector in the direction $e_j$. On these $|I|+(s-|I|)=s$ vectors, the matrix of the relative
differentials has the block form
\[
 \begin{pmatrix}
  \mathrm{Id}_{|I|} & 0\\
  * & \operatorname{diag}(z_j(x))_{j\notin I}
 \end{pmatrix}.
\]
Every diagonal entry in the second block is nonzero. The matrix is therefore invertible, so the relative differentials of all $s$ equations are independent. The ambient space $\mathcal T$ is smooth over $B$ of relative dimension $n+1+s$. The relative Jacobian criterion now shows that $W$ is smooth of relative dimension $n+1$ near $\widehat s(b)$. This holds for every $b$.

\noindent\emph{Rank and degree.}
Put $V=\widehat s^*T_{\mathcal T/B}$ and $E=\widehat s^*T_{W/B}$. The relative tangent sequence of the $(\mathbb{G}_m)^s$-torsor $\mathcal{T}\to\mathcal{Z}^\times$ restricting to $\widehat{s}$ gives the first of the following exact sequences; the differential of \eqref{eq:frames}, whose surjectivity was just proved, gives the second:
\begin{align*}
 0&\longrightarrow\cO_B^{\oplus s}\longrightarrow V \longrightarrow E_0\longrightarrow0,\\
 0&\longrightarrow E\longrightarrow V \longrightarrow\bigoplus_{i=1}^s A_i\longrightarrow0.
\end{align*}
In particular $\rk (V)=n+1+s$ and $\rk (E)=n+1$. Taking determinants in these two sequences and then in \eqref{eq:cone-tangent} gives
\[
 \det E
 \simeq\det V\otimes\bigotimes_i A_i^{-1}
 \simeq\det E_0\otimes\bigotimes_i A_i^{-1}
 \simeq f^*\cO_X(-K_X-D).
\]
Consequently
\[
 \deg E=-K_X\cdot C-\sum_i D_i\cdot C
       =-(K_X+D)\cdot C=\delta.
\]

\noindent\emph{Forgetting the frames away from the contacts.}
Over $B\setminus f^{-1}D$, every $b_i$ is a nonvanishing section. Equation \eqref{eq:frames} then implies that $p^*d_i$ is nonvanishing as well. Conversely, if all the $p^*d_i$ are nonvanishing, there is a unique isomorphism $u_i$ sending $p^*d_i$ to $q^*b_i$. This construction works over arbitrary test schemes and is compatible with base change. It gives the isomorphism of schemes as claimed.
\end{proof}

From now on, we replace $W$ by its relative smooth locus. This is an open neighborhood of $\widehat s(B)$; by the last assertion of Lemma \ref{lem:W}, it also contains the entire locus over $B\setminus f^{-1}D$. The family is separated and of finite type, being an open subscheme of the
zero scheme \eqref{eq:frames} in a product of line-bundle frame torsors. It therefore satisfies all the assumptions of Section~\ref{sec:jets}, with relative dimension $r=n+1$ and degree of the relative tangent bundle being $\delta$.

\section{From finite jets to exact polynomial equations}\label{sec:polynomial}

The jet given by Lemma~\ref{lem:positive} is defined only modulo the $(k+1)$-st power of the fiber parameter. To obtain a map on an actual surface, we will prove that its coordinate polynomials have
such small degrees that every relevant truncated identity is already an exact identity. This applies both to the equations of $X$ and to the ideals of the divisor components.

\subsection{Choice of the jet order and its coefficient bundles}

Choose finite sets of homogeneous generators
\[
 F_\alpha\quad\text{for }I_X\subset\C[z_0,\ldots,z_N],\qquad
 G_{ij}\quad\text{for }I_{D_i}\subset\C[z_0,\ldots,z_N].
\]
Here $I_{D_i}$ is the homogeneous ideal of $D_i$ as a closed subscheme of $\PP^N$. The use of ideals, rather than merely set-theoretic equations, will preserve the full contact multiplicities.
Let 
\begin{equation}\label{eqn:deg_poly_bound}
\Delta\geq1
\end{equation}
be an upper bound for the degrees of all these generators. Since the harmonic numbers satisfy $h_k\to\infty$, we can choose $k\geq2$ such that
\begin{equation}\label{eq:choosek}
 \delta h_k>2r\Delta a.
\end{equation}
An explicit permissible choice is
\begin{equation}\label{eq:explicit-k}
k=\max\{2,\lceil \exp(c_0)\rceil\}, \quad \text{where } c_0=\frac{2r\Delta a}{\delta}.
\end{equation}
This is because $h_k\geq\log(k+1)>c_0$.

Applying Lemma~\ref{lem:positive} to the family constructed in Section~\ref{sec:frames}. Let $\rho:B'\to B$ and $L$ be its finite cover and positive parameter line bundle. Write
\begin{equation}\label{eqn:degrees2}
 d=\deg\rho,\qquad \mu=\frac{\deg L}{d},\qquad A'=\rho^*A.
\end{equation}
The coordinate functions of the cone, restricted to the jet, have unique expressions
\begin{equation}\label{eq:P}
 P_j(t)=\rho^*f_j+\sum_{q=1}^k P_{jq}t^q, \qquad P_{jq}\in H^0(B',A'\otimes L^{-q}).
\end{equation}
To see both the uniqueness and the indicated coefficient bundles, let
\[
 S_k=\Spec_{B'}\left(\bigoplus_{q=0}^k L^{-q}\right)
\]
with truncated multiplication. A section of the pullback of $A'$ to $S_k$ is, by the affine projection formula, an element of
\[
 H^0\left(B',A'\otimes\bigoplus_{q=0}^k L^{-q}\right).
\]
Its summands are precisely the coefficients in \eqref{eq:P}.

For a coordinate description, choose local sections $\alpha$ of $A'$ and $\varepsilon$ of $L$. Write a fiber vector of $L$ as $t\varepsilon$, and write $P_{jq}=p_{jq}\alpha\otimes\varepsilon^{-q}$. On changing the sections to $\alpha'=v\alpha$ and $\varepsilon'=u\varepsilon$, we have
\[
t'=u^{-1}t,\quad p'_{jq}=v^{-1}u^q p_{jq},\quad \alpha'\sum_qp'_{jq}(t')^q=\alpha\sum_qp_{jq}t^q.
\]
Thus the finite sums in \eqref{eq:P} also define global polynomial sections on $\Tot(L)$, not just expressions in unrelated local parameters. The issue is whether these polynomial sections satisfy
the equations beyond order $k$.

We now work as in \cite[Lemma 4.1]{GX}. In the notations of (\ref{eqn:some_numbers}) and (\ref{eqn:degrees2}), we have 
\[
 \deg(A'\otimes L^{-q})=da-q\deg L=d(a-q\mu).
\]
A line bundle of negative degree on a smooth projective curve has no nonzero sections: the divisor of a nonzero section would be effective and would have that negative degree. Therefore
\[
 P_{jq}=0\in H^0(B',A'\otimes L^{-q}) \quad \text{whenever } \deg(A'\otimes L^{-q})<0, \text{ i.e. } q>\frac{a}{\mu}.
\]
Set 
\begin{equation}\label{eqn:Nk}
N_k=\lfloor a/\mu\rfloor.
\end{equation}
By \eqref{eq:slope} and the choice \eqref{eq:choosek}, we have
\begin{equation}\label{eq:degreebound}
 \frac{a}{\mu}<\frac{2rka}{\delta h_k}<\frac{k}{\Delta},  \quad\text{hence}\qquad \Delta N_k<k.
\end{equation}
Every $P_j(t)$ thus has fiber degree at most $N_k$. The strict inequality in \eqref{eq:degreebound} is the reason that a finite jet will suffice for all the equations chosen above.

\subsection{Exact equations of the target}

Suppose the degree of $F_\alpha$ is $e_\alpha$. Substituting \eqref{eq:P} into $F_\alpha$ gives a polynomial section of $\pi_L^*(A')^{e_\alpha}$, where 
$$\pi_L:\Tot(L)\to B'$$ 
is the projection. Its fiber degree is at most 
$$e_\alpha N_k\leq\Delta N_k<k.$$ 
In local line-bundle sections over an affine open $U=\Spec R\subset B'$, write it as
\[
 F_\alpha(P(t))=\sum_{q=0}^{e_\alpha N_k}c_qt^q.
\]
The original jet maps to the cone over $X$. Consequently, the image of this polynomial in $R[t]/(t^{k+1})$ is zero. The elements $1,t,\ldots,t^k$ form a free $R$-basis of that quotient. All the coefficients $c_q$ are therefore $0$, and $F_\alpha(P(t))=0$ as an actual polynomial. This reasoning applies to every generator and on every such open $U$.

It follows that $P(t)$ gives a morphism from $\Tot(L)$ to the twisted affine cone. This is the polynomial realization of \cite[Theorem 4.2]{GX}; the ideals for divisor components will be treated by the same coefficient argument below. Where the tuple $P(t)$ is nonzero, it can be projectivized, giving a rational map
\begin{equation}\label{eq:Phi}
 \Phi:\Tot(L)\dashrightarrow X.
\end{equation}
Its domain contains the entire zero section $\sigma$. Indeed, $P_j(0)=\rho^*f_j$, and these constant terms have no simultaneous zero. Thus the open complement of the simultaneous zero scheme of the $P_j$ contains $\sigma(B')$, and on that section $\Phi=f\rho$.

\subsection{The boundary ideal and its full contact multiplicity}
The following Lemma describes the behavior of the divisor components near the zero section of $\text{Tot}(L)$.

\begin{lemma}[Exact vertical divisors near the section]\label{lem:boundary}
Put $T_i=(f\rho)^*D_i$ on $B'$. In a neighborhood of $\sigma$ where $\Phi$ is regular, we have, after possibly shrinking that neighborhood,
\begin{equation}\label{eq:vertical}
\Phi^*D_i=\pi_L^*T_i\qquad\text{for every }i,
\end{equation}
where $\pi_L:\Tot(L)\to B'$ is the projection.
\end{lemma}

\begin{proof}
Fix $i$. We first prove divisibility by the complete local equation of $T_i$, and then prove that the remaining factor is a unit near the zero section. Both points are needed for equality of Cartier divisors.

\noindent\emph{Divisibility on the truncated source.}
Fix $y\in B'$ and choose an affine neighborhood $U=\Spec R$ on which the relevant line bundles are trivial. After shrinking $U$, we may also choose a homogeneous coordinate $z_\ell$ such that
$\rho^*f_\ell$ is nonvanishing on $U$. Trivializing $M_i$ near $f\rho(y)$ and writing its canonical section as $z_i$. In the induced sections of $(f\rho)^*M_i$ on $U$, we write
\[
 g_i=z_i\circ f\rho.
\]
This is a local equation of $T_i$; it is a nonzero element of the domain $R$. By further shrinking $U$, the reduction of the entire truncated jet over $U$ lies in the chosen target neighborhood. A
nilpotent thickening has the same underlying topological space as its reduction, so the truncated jet itself lies in that neighborhood.

On this jet, the frame equation reads
\[
 e_i(t)z_i(t)=g_i \quad\text{in }R[t]/(t^{k+1}),\qquad
 e_i(t)\in\bigl(R[t]/(t^{k+1})\bigr)^\times.
\]
Thus, its pulled-back ideal is exactly $(g_i)$ in this truncated ring. For example, if a local parameter $u$ at a contact point gives $g_i=u^{m_i}v_i$ with $v_i$ a unit, the identity is
\[
 z_i(t)=u^{m_i}v_i e_i(t)^{-1}\pmod{t^{k+1}}.
\]
In particular, the factor $u^{m_i}$ is kept.

Let $e_{ij}=\deg G_{ij}$. On the chart $z_\ell\ne0$, the function $G_{ij}/z_\ell^{e_{ij}}$ belongs to the ideal of $D_i$. Its pullback to the truncated jet is therefore divisible by $g_i$. The element $P_\ell(t)$ is invertible in $R[t]/(t^{k+1})$, since its constant term is a unit. Multiplication by its $e_{ij}$-th power gives
\begin{equation}\label{eq:divmod}
 G_{ij}(P(t))\in(g_i)\pmod{t^{k+1}}.
\end{equation}

\noindent\emph{Divisibility as an actual polynomial.}
The polynomial on the left of \eqref{eq:divmod} has degree at most 
$$e_{ij}N_k\leq\Delta N_k<k.$$ 
Its image in $(R/(g_i))[t]/(t^{k+1})$ is $0$. The basis $1,t,\ldots,t^k$ over $R/(g_i)$ shows that every coefficient lies in $(g_i)\subset R$. Dividing each coefficient by $g_i$ gives an actual polynomial $Q_{ij}(t)\in R[t]$ such that
\begin{equation}\label{eq:divexact}
G_{ij}(P(t))=g_i Q_{ij}(t).
\end{equation}
This is an identity in $R[t]$, with no truncation.

\noindent\emph{A generator whose quotient is a unit.}
Suppose first that $x=f\rho(y)$ belongs to $D_i$. The restrictions of the functions $G_{ij}/z_\ell^{e_{ij}}$ generate the ideal $\mathcal I_{D_i,x}$ in $\cO_{X,x}$. This ideal is an invertible module of rank one. Its quotient by the maximal ideal of $\cO_{X,x}$ has dimension one over the residue field. At least one of the generators considered here has a nonzero image in that quotient; Nakayama's lemma implies that this generator, say $G_{ij}/z_\ell^e$, generates $\mathcal I_{D_i,x}$ by itself. Equivalently, after shrinking the target neighborhood,
\[
 \frac{G_{ij}}{z_\ell^e}=z_i h, \quad h\in\cO_X^\times.
\]

We now work in the local domain $R_y=\cO_{B',y}$. Restricting this identity to the zero section gives
\[
 \frac{G_{ij}(P(0))}{P_\ell(0)^e}=g_i v, \quad v=(f\rho)^*h\in R_y^\times.
\]
On the other hand, \eqref{eq:divexact} identifies its left-hand side with 
$g_i Q_{ij}(0)/P_\ell(0)^e$. The element $g_i$ is nonzero in $R_y$, so it may be cancelled. We thus obtain
\[
 \frac{Q_{ij}(0)}{P_\ell(0)^e}=v\in R_y^\times.
\]
The regular function $Q_{ij}(t)/P_\ell(t)^e$ is therefore a unit in some neighborhood of $\sigma(y)$. In that neighborhood, the pullback of our chosen local generator of $D_i$ is
\[
 \Phi^*\left(\frac{G_{ij}}{z_\ell^e}\right)=g_i\frac{Q_{ij}(t)}{P_\ell(t)^e}.
\]
Its second factor is a unit, so its Cartier divisor is exactly the
pullback of $T_i$. This proves \eqref{eq:vertical} near
$\sigma(y)$, including its multiplicity.

If $f\rho(y)\notin D_i$, then both divisors in \eqref{eq:vertical} vanish near $\sigma(y)$: first shrink the base away from $\Supp T_i$, and then shrink the domain of $\Phi$ away from $\Phi^{-1}D_i$. Taking the union of the neighborhoods constructed for all $y\in B'$ gives an open neighborhood of the whole zero section on which \eqref{eq:vertical} holds for this $i$. Finally, intersect these neighborhoods over the finitely many indices $i$. The resulting open neighborhood still contains every point of $\sigma$, including all contact points, and all the equalities hold there simultaneously.
\end{proof}



The following lemma is the source of the rational curve that we eventually construct.

\begin{lemma}[Nonconstant generic parameter fiber]\label{lem:nonconstant}
The restriction of $\Phi$ to the generic fiber of $\Tot(L)\to B'$ is nonconstant.
\end{lemma}

\begin{proof}
Let $K=\C(B')$. Suppose that the rational map $\Phi$ on the generic fiber of $\text{Tot}(L)\to B'$ is constant. It is defined at $t=0$, so its constant value must be $f\rho$ at the generic point of $B'$. In local coordinate over $K$, equality of image points in $\mathbb{P}^N$ for $t=0$ and a nonzero $t$ gives the identities of polynomials
\[
 P_j(t)\rho^*f_\ell-P_\ell(t)\rho^*f_j=0 \quad\text{in }K[t],\qquad 0\leq j,\ell\leq N, \quad \text{where } f_0,...,f_N \text{ are in } (\ref{eqn:sections_of_A}).
\]
Taking the coefficient of $t^q$, for $q\geq1$, yields
\[
 P_{jq}\rho^*f_\ell=P_{\ell q}\rho^*f_j.
\]
These are equalities of global sections of $(A')^2\otimes L^{-q}$ on $B'$, since they hold at its generic point and $B'$ is integral.

The opens $U_\ell=\{\rho^*f_\ell\ne0\}$ cover $B'$. On $U_\ell$, define
\[
 \lambda_q=\frac{P_{\ell q}}{\rho^*f_\ell}.
\]
This is a regular section of $L^{-q}$. The cross-multiplication identities show that the definitions agree on all overlaps and that $P_{jq}=\lambda_q\rho^*f_j$ for every $j$. We have therefore constructed
\[
 \lambda_q\in H^0(B',L^{-q}),\quad  P_{jq}=\lambda_q\rho^*f_j.
\]
Since $\deg L>0$, the degree of $L^{-q}$ is negative for $q\geq1$. Thus every $\lambda_q$ vanishes. All positive-order coefficients of the cone jet vanish, so that jet is the constant section $s\rho$.

Let $U=B'\setminus\bigcup_i\Supp T_i$. It is a nonempty dense open subset. Over $U$, the last assertion of Lemma~\ref{lem:W} identifies $W$ with the divisor-complement part of the cone. The
frames are uniquely determined there, over the truncated source as well as over its reduction. Hence the original jet of $W$ is constant over $U$. At every point of $U$ all its positive-order
based-jet coefficients would be zero. This contradicts Lemma~\ref{lem:positive}, which gives a weighted representative whose coefficients do not all vanish at any point. This proves the Lemma.
\end{proof}

\section{The ruled-surface argument}\label{sec:ruled}

The polynomial map of Section~\ref{sec:polynomial} is regular near the whole zero section of $\Tot(L)$. In a neighborhood of the zero section, the pulled-back of the divisor $D$ is vertical, in the sense that it does not dominant the zero section. The role of the ruled surface is to turn this local statement along the zero section into a statement on a complete general parameter fiber. In particular, we take care of the vertical components by further shrinking $B'$.

The
positivity $\deg L>0$ is decisive: the only horizontal curve which can avoid the zero section is the infinity section.

\subsection{The two sections and their intersection numbers}

We use the convention that $\PP_{\mathrm{lines}}(V)$ parameterizes $1$-dimensional subspaces of the fibers of $V$. Put
\[
 S_0=\PP_{\mathrm{lines}}(\cO_{B'}\oplus L), \quad \pi_0:S_0\longrightarrow B'.
\]
The line $\C\oplus0$ defines the zero section $\sigma$, and the line $0\oplus L_y$ defines the infinity section $\tau$. To see the compactification explicitly, trivialize $L$ on an open subset $U\subset B'$ and write a line in $\cO\oplus L$ as $[a:b]$. On the chart $a\ne0$ the coordinate $t=b/a$ is a vector in $L_y$; hence this chart is $\Tot(L)|_U$, and $\sigma$ is given by $t=0$. The complementary chart $b\ne0$ has coordinate $u=a/b$, with $\tau$ given by $u=0$. On the overlap $u=t^{-1}$. Thus
\[
 S_0\setminus\tau=\Tot(L),\qquad \sigma\cap\tau=\varnothing.
\]

Let $\mathcal U=\cO_{S_0}(-1)$ denote the tautological subline of $\pi_0^*(\cO_{B'}\oplus L)$. Projections onto the two summands give
\[
 \mathcal U\longrightarrow\cO_{S_0},
 \qquad
 \mathcal U\longrightarrow\pi_0^*L.
\]
The first map vanishes precisely on $\tau$, and the second precisely on $\sigma$, each to order one. Interpreting these maps as sections of line bundles gives
\[
 \cO_{S_0}(\tau)\simeq\mathcal U^*,\qquad
 \cO_{S_0}(\sigma)\simeq\mathcal U^*\otimes\pi_0^*L.
\]
Consequently,
\begin{equation}\label{eq:section-line-bundle}
 \cO_{S_0}(\sigma-\tau)\simeq\pi_0^*L.
\end{equation}
This is the line-bundle identity used in \cite[proof of Proposition 6.1]{GX}.

Restriction of the tautological line gives $\sigma^*\mathcal U\simeq\cO_{B'}$ and
$\tau^*\mathcal U\simeq L$. Since the normal line bundle of an effective Cartier divisor $Z$ on a smooth surface is $\cO(Z)|_Z$, we obtain
\[
 N_{\sigma/S_0}\simeq L, \qquad N_{\tau/S_0}\simeq L^{-1}.
\]
Taking degrees, and taking first Chern classes in \eqref{eq:section-line-bundle}, we obtain
\begin{equation}\label{eq:sections}
 \sigma^2=\deg L>0,\qquad \tau^2=-\deg L,\qquad
 \sigma\cap\tau=\varnothing,
 \qquad [\sigma]=[\tau]+\pi_0^*c_1(L).
\end{equation}

\begin{lemma}\label{lem:onlyinfty}
An integral horizontal curve $Z\subset S_0$ disjoint from $\sigma$ is necessarily $\tau$.
\end{lemma}

\begin{proof}
Let $\nu_Z:\widetilde Z\to Z$ be the normalization. Horizontality here means that $\pi_0\nu_Z:\widetilde Z\to B'$ is nonconstant. It is therefore a finite surjective morphism of smooth projective curves. Write
\[
 d_Z=\deg(\pi_0\nu_Z)\geq1.
\]
The intersection with a line bundle is computed on the normalization, so
\[
 Z\cdot\pi_0^*c_1(L) =\deg\,((\pi_0\nu_Z)^*L) =d_Z\deg L>0.
\]
Suppose $Z\ne\tau$. Distinct integral curves on a smooth surface have nonnegative intersection number: their intersection is the sum of nonnegative local intersection multiplicities. Thus $Z\cdot\tau\geq0$. Using \eqref{eq:section-line-bundle}, we get
\[
 Z\cdot\sigma =Z\cdot\tau+d_Z\deg L>0.
\]
On the other hand, since $Z$ is disjoint from $\sigma$, we have $Z\cdot\sigma=0$. This contradiction proves $Z=\tau$.
\end{proof}

\subsection{Proof of Theorem~\ref{thm:main}}

We obtain the desired rational curve from the polynomial map in four steps.

\noindent\emph{Step 1: resolve away from the distinguished section.}
The rational map \eqref{eq:Phi} determines a rational map $S_0\dashrightarrow X$. Since $S_0$ is a smooth projective surface and $X$ is projective, resolving the indeterminacy of $\Phi$ by blowups of points  gives a diagram
\[
 S\xrightarrow{\ \beta\ }S_0, 
 \qquad h:S\longrightarrow X,
 \qquad \pi_1:=\pi_0\beta:S\longrightarrow B'.
\]
Here $h$ agrees with $\Phi\beta$ wherever it is defined. By Lemma~\ref{lem:boundary}, there is an open neighborhood $V$ of the whole zero section $\sigma\subset S_0$ on which $\Phi$ is regular and the identities \eqref{eq:vertical} hold. Therefore, we can choose the resolution to be an isomorphism over $V$: all the blowup centers lie outside its successive inverse images. We identify $\sigma$ with its strict transform in $S$. In particular,
\[
 \beta^{-1}(V)\simeq V, \qquad h|_\sigma=f\rho.
\]
Every exceptional curve of $\beta$ is vertical for $\pi_1$. Indeed, a point blowup creates a curve whose image in $B'$ is the image of its center; further blowups do not change the fact that each exceptional curve maps to a single point of $B'$.

\noindent\emph{Step 2: locate the horizontal boundary.}
The image of $h$ is not contained in any $D_i$, because its restriction to $\sigma$ is $f\rho$, whose image is $C\not\subset D$. Thus $h^*D$ is a well-defined effective Cartier divisor on the smooth integral surface $S$. Put $$T=\sum_i(f\rho)^*D_i.$$
On $\beta^{-1}(V)$, Lemma \ref{lem:boundary} gives an equality of Cartier divisors
\[
 h^*D=\pi_1^*T.
\]
We check explicitly that this excludes a horizontal component of $h^*D$ through any point of $\sigma$. Let $p=\sigma(y)$ and choose a local parameter $v$ on $B'$ at $y$. If $m$ is the multiplicity of $T$ at $y$, then in the regular local ring $\cO_{S,p}$, the ideal of $D$ is
\[
 I_{h^*D,p}=(v^m).
\]
For $m=0$ the divisor $D$ is absent near $p$. For $m>0$, its only codimension-$1$ support near $p$ is the local fiber $v=0$. Indeed, if a height-one prime $\mathfrak p$ contains $(v^m)$, then $v\in\mathfrak p$; the ruling is smooth near $p$, so $(v)$ is already a height-one prime and $\mathfrak p=(v)$. Thus an integral divisor component through $p$ must be that vertical fiber
component, and cannot be horizontal.

Let $Z'$ be any horizontal irreducible component of $h^*D$. We have proved $Z'\cap\sigma=\emptyset$. Since every exceptional curve is vertical, $\beta(Z')$ is an integral horizontal curve $Z\subset S_0$. Moreover, $Z\cap\sigma=\emptyset$: if a point of $Z$ lay on $\sigma$, its unique inverse image under $\beta$ would lie on the same section, contradicting $Z'\cap\sigma=\emptyset$. Lemma~\ref{lem:onlyinfty} therefore gives
\[
 \beta(Z')=\tau.
\]
In particular, every horizontal boundary component is the strict transform of $\tau$. This claim allows arbitrary divisor multiplicity on that section.

\noindent\emph{Step 3: restrict to good fibers.}
Let $$E_{\mathrm{bad}}\subset B'$$ 
be the union of the following finite sets: the images of all blowup centers, and the images of all vertical components of $h^*D$. 
Write $$U=B'\setminus E_{\mathrm{bad}}.$$ 
Over $U$, the morphism $\beta$ is an isomorphism and the ruling is the original $\PP^1$-bundle. For
every $y\in U$ its fiber $F_y$ is therefore $\PP^1$, with distinguished points $0=\sigma(y)$ and $\infty=\tau(y)$, and Step 2 gives
\[
 \Supp\bigl((h|_{F_y})^*D\bigr)\subseteq\{\infty\}.
\]
The divisor on the left is defined because $F_y$ is not a component of $h^*D$; all such vertical components have been removed.

Set $\mathcal M=h^*\cO_X(1)$. The degree of $\mathcal M|_{F_y}$ is constant on $U$. One may see this from the flat projective family $\pi^{-1}(U)\to U$: a line bundle is flat over the base, its fiberwise Euler characteristic is locally constant, and Riemann--Roch on $\PP^1$ gives
\[
 \deg(\mathcal M|_{F_y})=\chi(F_y,\mathcal M|_{F_y})-1.
\]
The open curve $U$ is connected. Its common fiber degree is positive, because the generic restriction of $h$ is nonconstant by Lemma~\ref{lem:nonconstant}, and the pullback of an ample line bundle under a nonconstant morphism of complete curves has positive degree. It follows that every map
\[
 h_y:=h|_{F_y}:\PP^1\longrightarrow X, \quad y\in U,
\]
is nonconstant. It satisfies $h_y(0)=f\rho(y)$ and the pull-back of $D$ under $h_y$ is supported at $\infty$.

Give the image $R_y=h_y(\PP^1)$ its reduced structure. It is an integral complete curve, not contained in $D$. Write $\nu_y:\widetilde R_y\to R_y$ for its normalization, normality of
$\PP^1$ implies that $h_y$ factors through $\nu_y$: 
\[
 h_y:\PP^1\xrightarrow{\ q_y\ }\widetilde R_y
 \xrightarrow{\ \nu_y\ }R_y\hookrightarrow X.
\]
The morphism $q_y$ is finite and surjective, of a positive degree $e_y$. By Riemann--Hurwitz, we have
\[
 2\cdot 0-2=e_y\bigl(2g(\widetilde R_y)-2\bigr) +\deg\operatorname{Ram}(q_y),
 \qquad \deg\operatorname{Ram}(q_y)\geq0.
\]
Hence $g(\widetilde R_y)=0$, and $\widetilde R_y\simeq\PP^1$. Furthermore, if 
$A_y=\nu_y^{-1}(D\cap R_y)$ is the set of intersection points with $D$ on the normalization, then
\[
 q_y^{-1}(A_y)=h_y^{-1}(D)\subseteq\{\infty\}.
\]
Surjectivity of $q_y$ implies that each point of $A_y$ has a nonempty inverse image, and distinct points have disjoint inverse images. Consequently, $$\#A_y\leq1.$$ 
This proves the required condition on the number of intersection points with $D$. 

\smallskip\noindent\emph{Step 4: descend the good locus to $C$.}
The map $f\rho:B'\to C$ is finite and surjective. The image of $E_{\mathrm{bad}}$ is therefore a finite closed subset of $C$. Also $C\cap D$ is finite, since $C\not\subset D$. Thus
\[
 C^\circ=C\setminus\bigl((C\cap D)\cup(f\rho)(E_{\mathrm{bad}})\bigr)\subset C\setminus D
\]
is open and dense. If $x\in C^\circ$, surjectivity gives a point $y\in B'$ with $f\rho(y)=x$; by the definition of $C^\circ$, such a point lies in $U$. The curve $R_y$ constructed in Step 3 passes through $x$ and has all the required properties. This completes the proof.

\subsection{The degree of the resulting curve}
Here we derive an estimate of the degree of the rational curve constructed in the proof of Theorem \ref{thm:main}.

In the above proof of Theorem~\ref{thm:main}, we fix the embedding $X\hookrightarrow\PP^N$ in (\ref{eqn:embedding}) and put $a=\deg f^*\cO_X(1)$, $r=\dim X+1$, and $\delta=-(K_X+D)\cdot C$ as in (\ref{eqn:some_numbers}). Let $\Delta\geq1$ be an upper bound of the degrees of the chosen homogeneous generators of $I_X$ and every $I_{D_i}$, as in (\ref{eqn:deg_poly_bound}).

\begin{proposition}[Explicit ample-degree bound]\label{prop:ample-bound}
For any integer $k\geq2$ satisfying \eqref{eq:choosek}, the curves in Theorem \ref{thm:main} can be chosen, through every point of a dense open subset of $C$, with
\begin{equation}\label{eq:ample-degree-bound}
 \deg_{\cO_X(1)}R_y
 \leq\deg h_y^*\cO_X(1)
 \leq\left\lfloor\frac a\mu\right\rfloor
 <\frac{2rka}{\delta h_k}.
\end{equation}
In particular, the choice \eqref{eq:explicit-k} gives the explicit integer bound
\begin{equation}\label{eq:integer-degree-bound}
 \deg_{\cO_X(1)}R_y
 \leq\left\lceil\frac{2rka}{\delta h_k}\right\rceil-1.
\end{equation}
Here $\deg_{\cO_X(1)}R_y$ means $c_1(\cO_X(1))\cdot R_y$ for the reduced integral image curve, and $\mu=\deg L/\deg\rho$ as in (\ref{eq:slope}).
\end{proposition}

\begin{proof}
We apply the degree argument of \cite[Corollary 4.3]{GX}. Fix a fiber $F_y$ in the proof of Theorem \ref{thm:main} and trivialize the $1$-dimensional vector spaces $L_y$ and $A_{\rho(y)}$. The coordinate sections then become polynomials
\[
 P_0(t),\ldots,P_N(t)\in\C[t],
 \qquad \deg P_j\leq N_k:=\left\lfloor\frac a\mu\right\rfloor.
\]
Their constant terms are not all zero. On $\PP^1$ with homogeneous coordinates $[S:T]$ and affine coordinate $t=T/S$, we homogenize them to a common degree:
\[
 H_j(S,T)=S^{N_k}P_j(T/S)\in\C[S,T]_{N_k}.
\]
Let $G(S,T)$ be a greatest common homogeneous divisor of the nonzero $H_j$, and write
\[
 H_j=G\,\overline H_j,\qquad  m=N_k-\deg G\leq N_k.
\]
This division removes the common zeros on the whole of $\PP^1$, including any common zero at infinity. More precisely, over $\C$ a common zero $[S:T]$ of the $\overline H_j$ would give a common linear factor, contrary to the definition of $G$. Hence the $\overline H_j$ have no simultaneous zero. They define a morphism
\[
 \overline h_y:\PP^1\longrightarrow\PP^N,
 \qquad [S:T]\longmapsto [\overline H_0(S,T):\cdots:\overline H_N(S,T)].
\]
The associated basepoint-free quotient $\cO_{\PP^1}^{\oplus(N+1)}\twoheadrightarrow\cO_{\PP^1}(m)$
identifies
\[
 \overline h_y^*\cO_{\PP^N}(1)\simeq\cO_{\PP^1}(m).
\]
On the dense open subset where the original polynomial tuple defines a point in $\mathbb{P}^N$, $\overline h_y$ agrees with $h_y$. Since $\PP^N$ is separated, the two morphisms agree everywhere. In particular, $\overline h_y$ has image in $X$, and
\[
 \deg h_y^*\cO_X(1)=m\leq N_k.
\]

The degree of the parametrized curve must be distinguished from the degree of its reduced image. For the factorization through $\nu_y:\widetilde R_y\to R_y$ used above, the projection formula gives
\[
 \begin{aligned}
 \deg h_y^*\cO_X(1)
 &=\deg q_y^*\nu_y^*(\cO_X(1)|_{R_y})\\
 &=e_y\deg\nu_y^*(\cO_X(1)|_{R_y})\\
 &=e_y\,\bigl(c_1(\cO_X(1))\cdot R_y\bigr),
 \qquad e_y=\deg q_y\geq1.
 \end{aligned}
\]
Thus $\deg_{\cO_X(1)}R_y\leq\deg h_y^*\cO_X(1)$. The strict upper bound in (\ref{eq:ample-degree-bound}) follows from \eqref{eq:slope}:
\[
 N_k\leq\frac a\mu<\frac{2rka}{\delta h_k}.
\]
An integer strictly less than a real number $T$ is at most $\lceil T\rceil-1$; this proves \eqref{eq:integer-degree-bound}. The bound holds for every fiber over $U$.

Finally, the estimate does not require a bound on the degree of the auxiliary cover. To check this directly, write $d=\deg\rho$. The coefficient of order $q$ in \eqref{eq:P} is a section of a line bundle of degree
\[
 \begin{aligned}
 \deg(\rho^*A\otimes L^{-q})
 &=d\,a-q\deg L=d\,(a-q\mu).
 \end{aligned}
\]
Since $d>0$, negativity is equivalent to $q>a/\mu$, independently of $d$. This is exactly the coefficient cutoff used to define $N_k$. \eqref{eq:slope} is likewise already normalized by $d$. Therefore, after fixing the embedding and its generators, the choice \eqref{eq:explicit-k} and the bound \eqref{eq:integer-degree-bound} depend only on $r,a,\delta,\Delta$.
\end{proof}

\section{Application to uniruledness}\label{sec:uniruled}

\begin{proof}[Proof of Proposition \ref{prop:logbb-quotient-surfaces}]
Choose a projective birational morphism
\[
\pi:X\longrightarrow S
\]
from a smooth projective surface such that
\[
B=\bigl(\pi^{-1}(D\cup\operatorname{Sing}(S)\cup V)\bigr)_{\mathrm{red}}
\]
is a simple normal crossings divisor. Such a morphism is obtained by resolving the surface and the boundary divisor and blowing up the remaining isolated points of the indicated closed subset. We choose $\pi$ to be an isomorphism outside $D\cup\operatorname{Sing}(S)\cup V$.

Since $S$ has quotient singularities, $K_S+D$ is $\mathbb Q$-Cartier. The nonexceptional components of $B$ are precisely the strict transforms of the components of $D$, with coefficient one. Hence
\begin{equation}\label{eq:logbb-resolution-formula}
    K_X+B=\pi^*(K_S+D)+\sum_i b_iE_i,
\end{equation}
where the $E_i$ are the exceptional prime divisors. We do not need to know the signs of the coefficients $b_i$.

Let $\Sigma\subset S$ be a finite set containing $\operatorname{Sing}(S)$, the singular points of $D$, $V$, and the images of all exceptional curves of $\pi$. For sufficiently large $m$, a general member
\[ 
C_0\in |mH|
\]
is smooth and irreducible, avoids $\Sigma$, and meets $D$ transversely in its smooth locus. We may also assume that $g(C_0)>0$, since
\[
2g(C_0)-2=m^2H^2+mK_S\cdot H
\]
is positive for sufficiently large $m$.

Write $C\subset X$ for the strict transform of $C_0$. Then $C\cong C_0$ and $C$ is disjoint from every exceptional curve. Intersecting \eqref{eq:logbb-resolution-formula} with $C$ gives
\begin{equation}\label{eq:logbb-negative-test-curve}
 (K_X+B)\cdot C =(K_S+D)\cdot C_0 =m(K_S+D)\cdot H <0.
\end{equation}
Notice that $C$ may meet $B$; only its disjointness from the exceptional locus is used in this calculation.

Applying Theorem \ref{thm:main} to $(X,B)$ and $C$, we obtain, through a general point of $C\setminus B$, an irreducible rational curve $L\subset X$ whose normalization map
\[
 f_L:\mathbb P^1\longrightarrow L\hookrightarrow X
\]
satisfies
\begin{equation}\label{eq:logbb-upstairs-contact}
\#\operatorname{Supp}(f_L^*B)\leq 1.
\end{equation}
In particular, $L$ is not contained in $B$.

We next verify that these curves contain a dominating algebraic family. Fix an ample divisor $A$ on $X$. For each $e>0$, consider maps $f:\PP^1\to X$ with $\text{deg}(f^*A)=e$. These maps are parametrized by the Hom scheme
\[
\mathcal H_e=\operatorname{Hom}_e(\mathbb P^1,X).
\]
Each $\mathcal H_e$ is of finite type. Let
\[
\mathcal H_e\xleftarrow{\ q\ }\mathbb P^1\times\mathcal H_e\xrightarrow{\ \operatorname{ev}\ } X,
\]
be the universal evaluation map (at $0\in \PP^1$) and the projection. The locus parametrizing maps $f:\PP^1\to X$ with $\text{deg}(f^*B)=0$ is
\[
 \mathcal H_e^\circ =\mathcal H_e\setminus q\bigl(\operatorname{ev}^{-1}(B)\bigr).
\]
This is because maps $f$ with $\text{deg}(f^*B)=0$ precisely those whose images are disjoint from $B$. $\mathcal H_e^\circ$ is open because $q$ is proper. In particular, it is of finite type.

For $d>0$, the locus parametrizing maps $f:\PP^1\to X$ with $\text{deg}(f^*B)=d$ is defined from $\mathcal{H}_e$ as follows. We first restrict to maps $f$ whose images are not contained in $B$ and for which $\deg(f^*B)=d$. After marking a point $p\in\mathbb P^1$, impose the condition
\[
f^*B=d[p].
\]
This defines an algebraic locus of finite type. 

Thus, as $e>0$ and $d\geq 0$ vary, maps satisfying \eqref{eq:logbb-upstairs-contact} are parameterized by countably many varieties of finite type.

Suppose that no irreducible family among these parameter spaces dominates $X$. The closure of the evaluation image of each family is then a single rational curve. Since $g(C)>0$, each such curve meets $C$ in finitely many points. Consequently, all these families together could account for only countably many points of $C(\mathbb C)$. This contradicts the existence of such rational curves  through points in a nonempty open subset of $C$. Hence one of these families dominates $X$.

A general member of this dominating family is not contracted by the log resolution $\pi$. Its image $R\subset S$ is rational. Moreover,
\[
 \pi^{-1} \bigl(D\cup\operatorname{Sing}(S)\cup V\bigr) =|B|
\]
set-theoretically. The inverse image of the indicated closed subset under the composite map
$\mathbb P^1\to X\to S$ therefore has at most one point. Factoring this composite through the normalization of $R$ cannot increase the number of inverse-image points. This proves \eqref{eq:logbb-one-contact} and establishes the first assertion.

For the second assertion, set
\[
W=\operatorname{Sing}(S)\cup(V\setminus D).
\]
The pair $(S,D)$ is purely log terminal at every point of $W$. This is assumed at the singular points of $S$; at the remaining points, $S$ is smooth and the divisor is absent.

The first assertion gives a dominating family of rational curves meeting
$D\cup\operatorname{Sing}(S)$ at most once. By \cite[Corollary~5.9.3]{KM}, this family can
be replaced by a dominating family disjoint from $W$ while preserving that contact bound.
The resulting curves lie entirely in $S^{\mathrm{sm}}$, avoid $V\setminus D$, and meet $D$ at most once. This proves \eqref{eq:logbb-smooth-locus}.
\end{proof}

{\large \part{Campana curves}}
In this part, we establish Theorem \ref{cp:thm:main}.

\section{An auxiliary compactification}\label{cp:sec:aux}

We construct in this section the auxiliary pair explicitly. 

\begin{proposition}\label{cp:prop:aux}
For the pair $(X, \De)$ where $\De=\sum_{i=1}^s\left(1-\frac1{m_i}\right)D_i$ is as in \eqref{cp:eq:pair}, there are a normal projective variety $Y$, a surjective morphism $f:Y\to X$, and a reduced divisor $B\subset Y$ such that:
\begin{enumerate}
\item $K_Y+B$ is $\Q$-Cartier and $K_Y+B\sim_{\Q}f^*(K_X+\De)$;

\item $U=Y\setminus B$ is smooth;

\item for every geometric point $x\in X\setminus D$, $U_x\simeq(\Gm)^s$;

\item there are integral reduced Weil divisors $E_i$ with $f^*D_i=m_iE_i$. In particular, $E_i|_U$ is Cartier.
\end{enumerate}
Here, a component $E_i$ with $m_i=1$ may be kept or omitted. If all such components are omitted, $s$ in \textup{(iii)} is replaced by the number of remaining components.
\end{proposition}

\begin{proof}
Choose Cartier divisors $A_i$ on $X$, and put
\[
 L_i=m_iA_i-D_i,\qquad \mathcal L_i=\OO_X(L_i).
\]
We can choose $A_i$ to be ample, but their positivity will play no role. Let
\[
 q_i:W_i=\PP_X(\mathcal L_i\oplus\OO_X)\longrightarrow X
\]
be the projective bundle of $1$-dimensional quotients. Write $S_{i,0}$ for the section given by the quotient to $\mathcal L_i$, and $S_{i,\infty}$ for the section given by the quotient to $\OO_X$. Let $\xi_i$ denote the tautological class, then
\begin{equation}\label{cp:eq:sections}
 S_{i,0}\sim\xi_i,\qquad
 S_{i,\infty}\sim\xi_i-q_i^*L_i,\qquad
 K_{W_i}+S_{i,0}+S_{i,\infty}\sim q_i^*K_X.
\end{equation}
Form the fiber product
\[
 W=W_1\times_X\cdots\times_XW_s,\qquad q:W\longrightarrow X.
\]
To simplify notations, we use the same symbols for the pulled-back sections on $W$, and set
$$S=\sum_i(S_{i,0}+S_{i,\infty}).$$ 
The relative canonical bundle formula gives
\begin{equation}\label{cp:eq:relcan}
 K_W+S\sim q^*K_X.
\end{equation}
For each $i$, define
\begin{equation}\label{cp:eq:branch}
 F_i=S_{i,0}+(m_i-1)S_{i,\infty}+q^*D_i,
 \qquad
 \mathcal M_i=\OO_W(q^*A_i+S_{i,\infty}).
\end{equation}
By \eqref{cp:eq:sections}, $F_i\sim m_i(q^*A_i+S_{i,\infty})$. Choose a section $\sigma_i$ of $\mathcal M_i^{\otimes m_i}$ with zero divisor $F_i$. The associated finite locally free cyclic-cover algebra is
\[
 \mathcal A_i=\bigoplus_{a=0}^{m_i-1}\mathcal M_i^{-a},
\]
where multiplication is determined by $\sigma_i$. Let
\begin{equation}\label{cp:eq:jointcover}
 Z=\Spec_W\left(\bigotimes_{i=1}^s\mathcal A_i\right),
 \qquad \pi:Y\to Z\to W
\end{equation}
be its normalization. Thus we normalize the joint cyclic-cover algebra, rather than take an unspecified component of a fiber product of normal covers.

We first check integrality. Let $K=\C(X)$. Over the generic point of $X$, choose independent fiber coordinates $t_i$ and rational trivializations of the line bundles. The generic equations are
\begin{equation}\label{cp:eq:kummerfield}
 z_i^{m_i}=d_it_i,\qquad d_i\in K^*.
\end{equation}
Their joint field is $K(z_1,\ldots,z_s)$, with $t_i=z_i^{m_i}/d_i$. Its degree over $K(t_1,\ldots,t_s)$ is $\prod_i m_i$. Hence the generic fiber of the algebra in \eqref{cp:eq:jointcover} is a field. Since that algebra is locally free over the integral variety $W$, it embeds in its generic fiber on each affine open. It is therefore a domain. Normalization is finite for varieties over $\C$, so $Y$ is integral, normal, and projective, and $\pi$ is finite
surjective. 

Put $f=q\circ\pi$. Let $B_{i,0}$ and $B_{i,\infty}$ be the reduced inverse images in $Y$ of the two
sections, and let $E_i$ be the reduced inverse image in $Y$ of $q^*D_i$. At the generic point of each of the three branch divisors in \eqref{cp:eq:branch}, the $i$-th extension has ramification index $m_i$: its branch multiplicity is respectively $1$, $m_i-1$, or $1$. The other cyclic extensions are \'etale there. For $m_i=1$, this says simply that the corresponding extension is the identity. Thus
\begin{equation}\label{cp:eq:pullbackdiv}
 \pi^*S_{i,0}=m_iB_{i,0},\qquad
 \pi^*S_{i,\infty}=m_iB_{i,\infty},\qquad
 f^*D_i=m_iE_i.
\end{equation}
There is no other ramification in codimension $1$. Set
$$B=\sum_i(B_{i,0}+B_{i,\infty});$$ 
the divisors in this sum have distinct supports. The characteristic-zero ramification formula, followed by
\eqref{cp:eq:pullbackdiv}, gives
\begin{equation}\label{cp:eq:mainidentity}
\begin{split}
 K_Y+B
 &\sim_{\Q}\pi^*K_W+\sum_i m_i(B_{i,0}+B_{i,\infty})+\sum_i(m_i-1)E_i\\
 &=\pi^*\left(K_W+S+q^*\De\right)\sim_{\Q}f^*(K_X+\De).
 \end{split}
\end{equation}
In particular, the left-hand side is $\Q$-Cartier. This calculation takes place on normal varieties and only uses codimension-$1$ ramification; it does not assume that $Y$ is smooth.

It remains to check smoothness off $B$. Locally on $X$, trivialize the line bundles and write $x_i$ for an equation of $D_i$. On $W\setminus S$, the fiber coordinates $t_i$ are invertible, and the joint cover has equations
\begin{equation}\label{cp:eq:localeq}
 z_i^{m_i}=t_ix_i,\qquad t_i\in\Gm.
\end{equation}
Units in the local equations can be absorbed into $t_i$. At a point of $X$, those $x_i$ which vanish are part of an \'etale coordinate system, by the SNC assumption. The corresponding rows of the Jacobian in \eqref{cp:eq:localeq} are independent in the $x_i$ directions, since $t_i$'s are units. For the remaining indices, the derivatives in the $t_i$ directions are units. The equations therefore define a smooth scheme. This scheme is already normal, so normalization changes nothing on this open subset. It is precisely $U=Y\setminus B$.

If $x\notin D$, each $x_i(x)$ is nonzero. In the fiber over $x$, we may eliminate $t_i=z_i^{m_i}/x_i(x)$, with all $z_i$ invertible. This identifies $U_x$ with $(\Gm)^s$. At the generic point of $D_i$, the reduced inverse image in $U$ has coordinates $z_i=0$, with $t_i$ and the $z_j$, $j\ne i$, invertible; it is an integral torus over $\C(D_i)$. Every component of $E_i$ dominates $q^*D_i$ under the finite map $\pi$, and hence meets this open subset. Thus $E_i$ is integral. Finally, a Weil divisor on the smooth variety $U$ is Cartier; this proves the last assertion of the proposition.
\end{proof}

\begin{remark}\label{cp:rem:separate}
Separate cyclic covers are essential in the smoothness above. A single cover branched along a reducible divisor may have local equation $z^m=tx_1x_2$; its open part need not be smooth. The independent equations \eqref{cp:eq:localeq} avoid this difficulty.
\end{remark}

\begin{remark}\label{cp:rem:noYfano}
The proposition does not assert that $Y$ is log Fano, which is not necessary for us. The positivity needed below comes entirely from $-(K_X+\De)$ and the identity \eqref{cp:eq:mainidentity}. The construction is a variant of \cite[Section~3]{LT}, suited to the reduced section divisor.
\end{remark}

\section{From pointwise existence to a covering family}\label{cp:sec:families}

In this section, we give the parameter-space argument needed to use Theorem~\ref{thm:main}. 

\begin{lemma}\label{cp:lem:param}
Let $V$ be smooth projective and let $\Gamma$ be an effective Cartier divisor. Nonconstant maps $r:\PP^1\to V$ not contained in $\Gamma$ for which $\#\Supp(r^*\Gamma)\leq1$ are parameterized, with a choice of a intersection point with $\Gamma$, by countably many finite-type schemes over $\C$.
\end{lemma}

\begin{proof}
Fix an ample divisor $H_V$. For every integer $d\geq1$, the scheme $\Hom_d(\PP^1,V)$ of maps of $H_V$-degree $d$ is of finite type. We restrict to the open subset of maps not contained in $\Gamma$. Pulling back $\Gamma$ via the universal morphism gives a relative effective Cartier divisor on the trivial $\PP^1$-bundle over this open subset. It is finite and flat over the parameter scheme. The degree $e\geq0$ of this divisor is locally constant. On a fixed-degree locus $H_{d,e}$, it determines a morphism
\[
 \rho_{d,e}:H_{d,e}\longrightarrow \operatorname{Hilb}^e(\PP^1)=\Sym^e(\PP^1).
\]
Define $\delta_e:\PP^1\to\Sym^e(\PP^1)$ by sending $q$ to $e[q]$, with the evident interpretation for $e=0$. The fiber product
\begin{equation}\label{cp:eq:incidence}
 M_{d,e}=H_{d,e}\times_{\Sym^e(\PP^1)}\PP^1
\end{equation}
parameterizes pairs $(r,q)$ with $r^*\Gamma=e[q]$. There are countably many pairs $(d,e)$, and each $M_{d,e}$ has finitely many irreducible components.
\end{proof}

\begin{lemma}[Covering family]\label{cp:lem:covering}
Let $(V,\Gamma)$ be a smooth projective SNC pair. Suppose there is a smooth proper family $\mathcal C\to T$ of geometrically integral curves, with $T$ an integral finite-type $\C$-scheme, and a morphism $\mathcal C\to V$ such that:
\begin{enumerate}
\item the evaluation $\mathcal C\to V$ is dominant;

\item every fiber maps isomorphically to a curve $C_t\subset V$ not contained in $\Gamma$;

\item $(K_V+\Gamma)\cdot C_t<0$ for all $t\in T$.
\end{enumerate}
Then a finite-type family of maps $r:\PP^1\to V$ satisfying $\#\Supp(r^*\Gamma)\leq1$ has dominant evaluation. Its evaluation away from the possible boundary point is dominant onto $V\setminus\Gamma$.
\end{lemma}

\begin{proof}
Take the reduced irreducible components of the spaces in Lemma~\ref{cp:lem:param}. If none of their evaluations is dominant, the closures of their images form a countable collection of proper closed subsets $Z_1,Z_2,\ldots\subsetneq V$.

For each $j$, the locus
\[
 T_j=\{t\in T \, | \, C_t\subset Z_j\}
\]
is closed and proper. Indeed, its complement is the image in $T$ of the open subset of $\mathcal C$ mapping to $V\setminus Z_j$; this image is open because $\mathcal C\to T$ is smooth. It is nonempty because the evaluation is dominant. An irreducible complex variety is not the union of countably many proper closed subsets. Indeed, its smooth locus has a nonempty analytic open subset which is a Baire space, while the intersection of any proper algebraic closed subset with this open subset is nowhere dense. Choose $t\notin\bigcup_jT_j$. Each intersection $C_t\cap Z_j$ is finite.

By Theorem~\ref{thm:main}, there is a nonempty open subset of $C_t$ covered by images of rational curves parameterized by the $M_{d,e}$. That open subset is thus contained in the countable set
$\bigcup_j(C_t\cap Z_j)$, which is impossible. Some evaluation is therefore dominant. Its restriction away from the marked divisor point is still dominant on $V\setminus\Gamma$: for $e>0$ the removed point maps into $\Gamma$, and for $e=0$ removing one point from each nonconstant curve does not change the closure of the total image.
\end{proof}

\begin{remark}\label{cp:rem:nobound}
Lemma~\ref{cp:lem:covering} selects one finite-type family from a countable collection, without a uniform bound on the degree. Its proof uses the uncountability of $\C$ and applies Theorem~\ref{thm:main} only to complex varieties.
\end{remark}

\section{Negative logarithmic degree}\label{cp:sec:negative}

We again consider the auxiliary pair of Proposition~\ref{cp:prop:aux}. Choose a projective log resolution
\[
 p:V\longrightarrow Y
\]
which is an isomorphism over the SNC locus of $(Y,B)$. In particular, it is an isomorphism over $U=Y\setminus B$. Define
\begin{equation}\label{cp:eq:Gamma}
 \Gamma=\bigl(p^{-1}(B)\cup\Exc(p)\bigr)_{\red}.
\end{equation}
This notation denotes a reduced union: an exceptional component appearing in the inverse image of $B$ is counted only once. Then $\Gamma$ is SNC and $V\setminus\Gamma\simeq U$.

The pair $(Y,B)$ is SNC at the generic point of every divisor. Its non-SNC locus, and the image of the exceptional locus of $p$, are contained in a closed subset $T_Y\subset Y$ of codimension at least two. There is a $\Q$-divisor $F$, supported on the exceptional locus, such that
\begin{equation}\label{cp:eq:discrepancy}
 K_V+\Gamma\sim_{\Q}p^*(K_Y+B)+F.
\end{equation}
We do not need the signs of the coefficients of $F$.

\begin{lemma}\label{cp:lem:negative}
If $\mathcal H$ in (\ref{cp:eq:L}) is big, there is a family of smooth complete-intersection curves 
$C_t\subset Y\setminus T_Y$ whose strict transforms on $V$ satisfy the assumptions of Lemma~\ref{cp:lem:covering}.
\end{lemma}

\begin{proof}
Let $n=\dim Y$ and choose a sufficiently ample Cartier divisor $H$ on $Y$. For $n\geq2$, the intersection of $n-1$ general members of $|H|$ is a smooth integral curve, avoids $T_Y$, and is not contained in $B$. These assertions follow from Bertini on the smooth locus and the dimension estimate $\dim T_Y-(n-1)\leq-1$. We restrict the parameter space of these curves to the nonempty
open subset where these assertions hold. The universal curves form a smooth proper family with dominant evaluation on $Y$. Their strict transforms, identified with the curves themselves, have dominant evaluation on $V$.

For any such curve $C$, \eqref{cp:eq:discrepancy} and avoidance of the exceptional locus give
\begin{equation}\label{cp:eq:negative}
 (K_V+\Gamma)\cdot\widetilde C =(K_Y+B)\cdot C =-f^*\mathcal H\cdot H^{n-1}<0.
\end{equation}
We check the strict inequality as follows. By bigness, $\mathcal H\sim_{\Q}A+G$ with $A$ ample and $G$ effective. Since $f$
is surjective, $f^*G$ is effective and
$(f^*G)\cdot H^{n-1}\geq0$. A suitable multiple of $A$ has a nonzero
effective representative; its inverse image is a nonzero effective
Cartier divisor on $Y$. Its intersection with $H^{n-1}$ is positive.
Thus $f^*\mathcal H\cdot H^{n-1}>0$.
Only this mixed intersection inequality is used; no bigness of $f^*\mathcal H$
on $Y$ is required.

If $n=1$, normality makes $Y$ smooth and the codimension-two locus is
empty. Take $C=Y$ as a constant family. The same degree calculation
applies. This proves the lemma in all dimensions under consideration.
\end{proof}

\Needspace{9\baselineskip}
\begin{proposition}\label{cp:prop:projected}
Suppose $\mathcal H$ is big. There is
an integral finite-type scheme $T$, a morphism
$h:\PP^1\times T\to X$, and a section $q:T\to\PP^1\times T$ such that:
\begin{enumerate}
\item each $h_t$ is nonconstant and its image is not contained in $D$;

\item evaluation of $h$ is dominant;

\item for every $i$, every contact of $h_t$ with $D_i$ away from $q(t)$ is divisible by $m_i$.
\end{enumerate}
\end{proposition}

\begin{proof}
Applying Lemmas~\ref{cp:lem:negative} and \ref{cp:lem:covering}, we obtain a dominant family $r:\PP^1\times T\to V$ with $r_t^*\Gamma=e[q(t)]$, for a fixed $e\geq0$. We then replace $T$ by a reduced irreducible component with dominant evaluation. Its evaluation on the complement of $q$ is dominant onto $U$.

Shrinking to the nonempty subset of $T$ corresponding to maps with images meeting $f^{-1}(X\setminus D)\cap U$. This subset is open because it is the image of an open subset under the smooth projection $\PP^1\times T\to T$. It is nonempty by dominance. 

Put $h=(f\circ p)\circ r$. Suppose $h_t$ were constant, with value $x$. The preceding restriction
implies that $x\notin D$. The nonconstant map $r_t$, restricted away from $q(t)$, would then give a nonconstant map $\A^1\to U_x\simeq(\Gm)^s$. Such a map does not exist: each coordinate would be a unit in $\C[t]$, hence a constant. The restriction really is nonconstant, since $r_t$ is nonconstant and $p$ is an isomorphism on $U$. Thus $h_t$ is nonconstant; it is not contained in $D$ by construction.

The evaluation of $h$ is dominant because the evaluation of $r$ is dominant and $f\circ p$ is surjective. On the complement of $q$, the map $p\circ r_t$ takes values in the smooth open subset $U$. There the identity $f^*D_i=m_iE_i$ is an identity of Cartier divisors. Pulling it back proves \textup{(iii)}.
\end{proof}

\section{Repairing contacts}\label{cp:sec:repair}

Here we obtain results on divisibility of contact orders after repairs.

\begin{lemma}[One-point case]\label{cp:lem:repair}
Let $h:\PP^1\to X$ be nonconstant and not contained in $D$. Suppose all contacts with $D_i$ away from a point $q$ are divisible by $m_i$. Set $N=\lcm(m_1,\ldots,m_s)$, with $N=1$ if $s=0$. Then thre is a map $\varphi:\PP^1\to\PP^1$ of degree $N$ such that $h\circ\varphi$ is a
divisible Campana rational curve. It has the same image as $h$.
\end{lemma}

\begin{proof}
Choose $b\notin h^{-1}(D)\cup\{q\}$ and an automorphism $a\in\PGL_2(\C)$ with $a(\infty)=q$ and $a(0)=b$. Let
\[
 \nu_N([u:v])=[u^N:v^N],\qquad \varphi=a\circ\nu_N.
\]
The map is totally ramified over $q$ and $b$ and is \'etale elsewhere. For every $p\in\PP^1$,
\begin{equation}\label{cp:eq:contactcover}
 \operatorname{ord}_{p}(h\circ\varphi)^*D_i
=e_p(\varphi)\,\operatorname{ord}_{\varphi(p)}h^*D_i.
\end{equation}
When $p=q$, the right-hand side of the above is a multiple of $N$, hence of $m_i$. It is $0$ when $p=b$. At every other point it is the original divisible contact. This proves the assertion, including when several $D_i$ pass through $h(q)$. Surjectivity of $\varphi$ proves the assertion about images.
\end{proof}

Lemma \ref{cp:lem:repair} concerns the map $h\circ\varphi$. It does not say that the normalization map of its image has divisible contact orders. Multiple covers are allowed in Definition~\ref{cp:def:campana}, as in \cite[Definition~2.2]{LT}. Since $\varphi$ has degree $N$, the numerical classes satisfy
\[
 (h\circ\varphi)_*[\PP^1]=N\,h_*[\PP^1].
\]
Thus the repair generally changes the class. The class kept fixed in Section~\ref{cp:sec:free} is the class of the repaired family.

\begin{lemma}[A dominant family]\label{cp:lem:familyrepair}
Let $(h,T,q)$ be as in Proposition~\ref{cp:prop:projected}. Then there is a finite-type family of divisible Campana rational curves whose evaluation dominants $X$.
\end{lemma}

\begin{proof}
View $q$ as a morphism $T\to\PP^1$. Introduce the frame space
\begin{equation}\label{cp:eq:frames}
 T'=\{(t,a)\in T\times\PGL_2(\C) \ |\ a(\infty)=q(t),\ h_t(a(0))\notin D\}.
\end{equation}
It is locally closed and of finite type. The projection $T'\to T$ is surjective: we can choose a second point outside the finite set $h_t^{-1}(D)\cup\{q(t)\}$ and then choose an automorphism with the two specified images. Define the family
\[
 G:\PP^1\times T'\longrightarrow X,\qquad G_{t,a}=h_t\circ a\circ\nu_N.
\]
It is an algebraic morphism. By Lemma~\ref{cp:lem:repair}, we have divisibility of contact orders fiberwise, and $G_{t,a}(\PP^1)=h_t(\PP^1)$. The union of the images is therefore unchanged, so evaluation map is dominant. We conclude by choosing a reduced irreducible component of $T'$ with dominant evaluation.
\end{proof}

\begin{lemma}[Fixing the contact type]\label{cp:lem:contacttype}
Let $G:\PP^1\times T\to X$ be a dominant family of nonconstant divisible Campana rational curves, with $T$ integral and of finite type. After shrinking $T$ and making a finite \'etale surjective base change, the union of the supports of the $G_t^*D_i$ is a disjoint union of marked sections, and
all their contact vectors are constant. The resulting family is a dominant family of nondegenerate stable log maps with fixed Campana contact data.
\end{lemma}

\begin{proof}
We first restrict $T$ to its smooth open dense subset. The images of the maps $G_t$ are not contained in any $D_i$. Hence each $G^*D_i$ is a relative effective Cartier divisor on 
$\PP^1\times T$, finite and flat over $T$. Its degree is constant after restricting to a
connected open subset. The reduced support of the union of these divisors is generically a finite \'etale scheme over $T$, since the function field has characteristic zero. Shrinking $T$ further so that this reduced support is finite \'etale over $T$. The finite \'etale scheme of orderings of geometric points of this reduced support is surjective over $T$. Choose a connected component of this cover; it remains finite \'etale and surjective, since $T$ is connected. On this component, the support splits into disjoint sections $p_1,\ldots,p_\ell$. Shrinking the original base once more if necessary. Then there are constant integers $c_{ji}\geq0$ such that
\begin{equation}\label{eqn:divisor_identity}
 G^*D_i=\sum_{j=1}^{\ell}c_{ji}p_j.
\end{equation}
Each contact vector $(c_{j1},\ldots,c_{js})$ is nonzero, and $m_i\, |\, c_{ji}$.

Give $X$ its divisorial log structure and the source its marked-curve log structure along the $p_j$. (\ref{eqn:divisor_identity}) induce the associated nondegenerate log map over the trivial log structure on the base. It is stable because the underlying maps are nonconstant maps from $\PP^1$; their automorphism groups are finite. The contact vectors are fixed by construction. Restriction to dense open subsets and surjective finite \'etale base change preserve dominance of evaluation. We may choose a component dominating the original base. If the support is empty, the same argument uses the trivial log structure on the source and no markings.
\end{proof}

\section{Freeness and rational connectedness}\label{cp:sec:free}

Here we prove that a nonempty open subfamily of the repaired family consists of free curves, with unchanged numerical class and contact vectors. We then apply the rational-connectedness theorem of
Chen--Lehmann--Tanimoto.

\subsection{The normal complex and the graph bundle}
Let $g:\mathbb{P}^1\simeq C\to X$ be nonconstant and not contained in $D$. Mark every point of $g^{-1}(D)$ and write $P\subset C$ for their reduced sum. Set
\[
 A=T_C(-P),\qquad E=g^*T_X(-\log D).
\]
The original logarithmic normal complex is
\[
 \mathcal N_g=[A\xrightarrow{\ d^{\log}g\ }E], \quad \text{in degrees } -1, 0.
\]
Since $g$ is nonconstant and the characteristic is zero, $d^{\log}g$ is generically nonzero and hence injective as a sheaf map. Its cokernel can have torsion at ramification points; it need not be a vector bundle.

Consider instead the graph
\[
 \widetilde g=(g,\operatorname{id}_C):C\longrightarrow (X,D)\times C,
\]
where the second factor has trivial log structure. Its logarithmic normal bundle $\mathcal N_{\widetilde g}$ is given by
\begin{equation}\label{cp:eq:graphnormal}
 0\longrightarrow A
 \xrightarrow{\ (d^{\log}g,\iota)\ }E\oplus T_C
 \longrightarrow\mathcal N_{\widetilde g}\longrightarrow0.
\end{equation}
Here $T_C$ is the ordinary tangent bundle and $\iota:A\to T_C$ is the natural inclusion. We check that the first map is a subbundle inclusion. Away from $P$, its second component is an isomorphism. Near a marked point (in $P$), choose a local parameter $s$ on $C$ and local equations $x_i$ for the divisor components through its image. Write
\[
 x_i\circ g=u_i(s)s^{c_i},\qquad u_i(0)\ne0.
\]
The logarithmic normal components of the image of $s\partial_s$ are
\[
 c_i+s\frac{u_i'(s)}{u_i(s)}.
\]
At least one $c_i$ is positive, hence nonzero over $\C$. The map is therefore nonzero on every tangent fiber, proving that its cokernel is locally free. This argument does not require $g$ itself to be an immersion.

Projection onto $E$ gives the distinguished triangle
\begin{equation}\label{cp:eq:normaltriangle}
 T_C\longrightarrow\mathcal N_{\widetilde g}
 \longrightarrow\mathcal N_g\longrightarrow T_C[1].
\end{equation}
Indeed, compare the complexes $[A\to E\oplus T_C]$ and $[A\to E]$, whose kernel is $T_C$ in degree $0$. Since $T_C\simeq\OO_C(2)$, for $a=1,2$ we have $H^1(C,T_C(-a))=H^2(C,T_C(-a))=0$. The cohomology exact sequence gives 
\begin{equation}\label{cp:eq:normalcomparison}
 H^1\bigl(C,\mathcal N_{\widetilde g}(-a)\bigr)
 \simeq\HH^1\bigl(C,\mathcal N_g(-a)\bigr),\qquad a=1,2.
\end{equation}
This is the comparison in \cite[Lemma~4.8, pp.~727--728, equation~(4.1)]{CLT}.

In the terminology of \cite[Definition~4.2, p.~726]{CLT}, $g$ is {\em free} if the group in \eqref{cp:eq:normalcomparison} vanishes for $a=1$, and {\em very free} if it vanishes for $a=2$. On $\PP^1$, the splitting theorem for vector bundles identifies these conditions with nefness and ampleness, respectively, of $\mathcal N_{\widetilde g}$. In particular, this definition is not an assertion that $g^*T_X(-\log D)$ is nef. If $P$ is empty, the same comparison applies without contact markings.

\subsection{An open subfamily of free curves}
We now construct a family of free curves in the above sense.

\begin{proposition}\label{cp:prop:fixedclassfree}
Let $T$ be a smooth integral variety, and let $G:\PP^1\times T\to X$ be a family of nonconstant maps with dominant evaluation map, none of whose images is contained in $D$. Suppose the preimages of components of $D$ are disjoint marked sections $P_1,\ldots,P_\ell$ and
\begin{equation}\label{eqn:divisor-identity2}
 G^*D_i=\sum_jc_{ji}P_j,
\end{equation}
where the $c_{ji}$ are constant and each contact vector is nonzero. Then there is a nonempty open subset $T_0\subset T$ such that every member over $T_0$ is free. The restricted family still has dominant evaluation map, and keeps the same numerical class and contact vectors of the original
family. In particular, divisible contact orders remain divisible.
\end{proposition}

\begin{proof}
We begin with some simplifying notations. Put 
$$S=\PP^1\times T \quad \text{and }P_S=\sum_jP_j.$$ For $t\in T$, write 
$$C=S_t,\quad g=G_t,\quad \text{and } P=P_S|_C.$$ 
The logarithmic tangent sequence, restricted to $C$, is
\begin{equation}\label{cp:eq:familytangent}
 0\longrightarrow T_C(-P)
 \longrightarrow T_S(-\log P_S)|_C
 \longrightarrow T_{T,t}\otimes_\C\OO_C\longrightarrow0.
\end{equation}
We explain this for completeness. Near a marked section, choose a relative coordinate $s$ with that section given by $s=0$. Logarithmic tangent vectors are generated by $s\partial_s$ and lifts of tangent vectors on $T$. This proves surjectivity on the right and identifies the kernel; away from
the markings it is the ordinary relative tangent sequence.

The divisor identity (\ref{eqn:divisor-identity2}) make
\[
 \widetilde G=(G,\operatorname{pr}_1): (S,P_S)\longrightarrow (X,D)\times\PP^1
\]
a logarithmic map. Its differential takes the left term of \eqref{cp:eq:familytangent} to the first map in \eqref{cp:eq:graphnormal}. Passing to quotients therefore gives
\begin{equation}\label{cp:eq:kappa}
 \kappa:T_{T,t}\otimes_\C\OO_C\longrightarrow\mathcal N_{\widetilde g}.
\end{equation}
Every parameter tangent vector thus gives a global section of the normal bundle of the graph of $g$. No global splitting of \eqref{cp:eq:familytangent} is needed: two local lifts differ by a section of $T_C(-P)$, whose image in $\mathcal{N}_{\widetilde g}$ vanishes.

At an unmarked point $u\in C$, identify the fiber of $\mathcal{N}_{\widetilde g}$ with $T_{X,g(u)}$ by
\begin{equation}\label{cp:eq:fiberidentification}
 [(v,w)]\longmapsto v-dg_u(w).
\end{equation}
Under this identification, the evaluation of \eqref{cp:eq:kappa} is the derivative of $G$ in the parameter directions. Indeed, a lift of a parameter tangent vector $v$ with source component $w$ has image $(dg_u(w)+d_TG(v),w)$; formula~\eqref{cp:eq:fiberidentification} cancels the first summand.

There are also global sections given by the sheaf map
\[
 T_C\longrightarrow\mathcal N_{\widetilde g},\qquad w\longmapsto[(0,w)].
\]
Their values under \eqref{cp:eq:fiberidentification} are $-dg_u(w(u))$. These sections come from the ordinary second target factor; the vector fields need not preserve the contact markings on the source. Since $T_{\PP^1}=\OO_{\PP^1}(2)$ is globally generated, these include all derivatives in the source-curve direction.

By dominance in characteristic zero, there is an unmarked point $(u,t)$ at which $dG$ is surjective. The parameter and source derivatives there span $T_{X,G(u,t)}$. Consequently, the sections just constructed defines a surjective map
\[
 H^0(C,\mathcal N_{\widetilde g})\longrightarrow\mathcal N_{\widetilde g}|_u.
\]
Write $\mathcal N_{\widetilde g}\simeq\bigoplus_j\OO_C(a_j)$. A negative summand has no global sections and would prevent this surjectivity. Hence every $a_j\geq0$, proving nefness and freeness.

Let $\Omega\subset S\setminus P_S$ be the nonempty open locus where $dG$ is surjective. Since the projection $\text{pr}_T:S=\mathbb{P}^1\times T\to T$ is smooth and hence open, $T_0=\operatorname{pr}_T(\Omega)$ is a nonempty open subset of $T$. Every fiber of $\text{pr}_T$ over $T_0$ meets $\Omega$, so the preceding argument proves that every such fiber is free. The source $\PP^1\times T_0$ is dense open in the integral variety $S$; restricting the dominant evaluation map to it preserves dominance.

Over $t\in T_0$, the map obtained is the original map $G_t$. The degree of the pullback of every line bundle on $X$ is constant in this connected proper flat family of curves; hence the numerical class is constant. The contact vectors are fixed by assumption. This proves the remaining assertions.
\end{proof}

\begin{remark}[The parameter-dimension convention]\label{cp:rem:slicing}
Proposition \ref{cp:prop:fixedclassfree} is the form of \cite[Proposition~4.4, p.~726]{CLT} needed here, with the additional conclusion about curve class. To compare our family with the
uniruledness convention in \cite[Definition~4.1, p.~726]{CLT}, put $n=\dim X$. That definition uses a base of dimension $n-1$.

At a general unmarked point $(u,t)$, the differential $dG$ is surjective and its restriction $dg_t$ to $T_{\PP^1,u}$ is nonzero. The latter follows from nonconstancy in characteristic zero. Choose an $(n-1)$-dimensional subspace $V\subset T_{T,t}$ whose image complements $dg_t(T_{\PP^1,u})$ in $T_{X,G(u,t)}$. Then
\[
 T_{\PP^1,u}\oplus V\longrightarrow T_{X,G(u,t)}
\]
is an isomorphism. Smooth local coordinates on $T$ give a smooth locally closed slice $W\subset T$ through $t$ with tangent space $V$: choose local functions whose independent differentials cut out $V$, and take their common zero locus near $t$. The restricted evaluation $\PP^1\times W\to X$ has surjective differential at $(u,t)$, so it is dominant and separable. When $n=1$, take $W=\{t\}$.

This restriction preserves the curve class and the contact data and satisfies \cite[Definition~4.1, p.~726]{CLT}. The same slicing argument occurs in \cite[Lemma~4.9, p.~729]{CLT}. Our direct proof of
Proposition~\ref{cp:prop:fixedclassfree} does not require the slice.
\end{remark}


\subsection{Proofs of the main statements}
We can now complete the proofs of results on Campana curves.

\begin{proof}[Proof of Theorem~\ref{cp:thm:main}]
Proposition~\ref{cp:prop:aux} constructs $(Y,B)$. Lemma~\ref{cp:lem:negative} and Lemma \ref{cp:lem:covering} give Proposition~\ref{cp:prop:projected}. Lemma~\ref{cp:lem:familyrepair} then
produces a covering family of divisible Campana rational curves. Lemma~\ref{cp:lem:contacttype} puts this family into fixed-contact form with a smooth integral base. Its numerical class is fixed.
Proposition~\ref{cp:prop:fixedclassfree} gives a nonempty open restriction on which every member is free and evaluation map remains dominant. Both the numerical class and divisibility are unchanged.
\end{proof}

A klt Campana orbifold is \emph{strongly Campana uniruled} if it has a free Campana rational curve whose numerical class lies in $\Int\Nef_1(X)$; this is \cite[Definition~5.9, p.~731]{CLT}. Here $\Nef_1(X)$ denotes the cone of numerical curve classes having nonnegative intersection with every effective divisor, and the interior is taken in $N_1(X)_{\mathbb R}$.

The following is \cite[Corollary~6.7]{CLT}.

\begin{theorem}[Chen--Lehmann--Tanimoto]\label{cp:thm:CLT}
Let $(X,\De)$ be a smooth projective SNC Campana orbifold over $\C$. If $X$ is rationally connected and $(X,\De)$ is strongly Campana uniruled, then $(X,\De)$ is Campana rationally connected.
\end{theorem}
In our application, strong Campana uniruledness will follow from the covering family and the Picard-number assumption.

\begin{lemma}\label{cp:lem:ordinaryRC}
For a smooth projective SNC weak Fano orbifold $(X,\De)$, the variety $X$ is rationally connected.
\end{lemma}

\begin{proof}
This result is standard; we include a proof for the reader's convenience.
Write $\mathcal H\sim_{\Q}A+G$ with $A$ ample and $G$ effective. Since $(X,\De)$ is klt, the pair $(X,\De+\delta G)$ remains klt for all sufficiently small $\delta\in\mathbb{Q}_{>0}$. Moreover,
\[
 -(K_X+\De+\delta G)
 \sim_{\Q}(1-\delta)\mathcal H+\delta A
\]
is ample when $0<\delta<1$. Thus $X$ is of Fano type. The rational-connectedness theorem for klt log Fano pairs \cite{HM} applies. Equivalently, it gives rational chain connectedness here, which is
rational connectedness for a smooth projective variety in characteristic zero. This ampleness calculation uses nefness of $\mathcal H$.
\end{proof}

\begin{proof}[Proof of Corollary~\ref{cp:cor:rhoone}]
First we recover the positivity needed for ordinary rational connectedness. Choose an ample Cartier divisor $H_X$. Since $\rho(X)=1$, write $\mathcal H\equiv aH_X$. A Kodaira decomposition $\mathcal H\sim_{\Q}A+G$, with $A$ ample and $G$ effective, gives $\mathcal H\cdot H_X^{\dim X-1}>0$, so $a>0$. Ampleness is numerical for $\Q$-Cartier divisors; hence $\mathcal H$ is ample and in particular nef. Thus Lemma~\ref{cp:lem:ordinaryRC} applies to $X$.

Let $\beta$ be the numerical class of the covering family in Theorem~\ref{cp:thm:main}. Nonconstancy implies $H\cdot\beta>0$ for any ample divisor $H$, so $\beta\ne0$. For every effective divisor $E$ on $X$, dominance of the evaluation map implies that there exists a member in the covering family not contained in $E$. The pullback of $E$ to this member is effective, whence $E\cdot\beta\geq0$. Thus $\beta\in\Nef_1(X)$.

When $\rho(X)=1$, a nonzero nef curve class lies in the interior of this $1$-dimensional cone. Proposition~\ref{cp:prop:fixedclassfree} gives a free Campana curve of class $\beta$. Hence $(X,\De)$ is strongly Campana uniruled. Lemma~\ref{cp:lem:ordinaryRC} gives ordinary rational connectedness of $X$, and Theorem~\ref{cp:thm:CLT} applies.
\end{proof}

\end{document}